\documentclass[preprint]{elsarticle}
           \usepackage[T2A]{fontenc}   
           \usepackage[utf8]{inputenc}
           \usepackage[english]{babel}
\usepackage{amsmath}
\usepackage{amsfonts}
\usepackage{amssymb}
\usepackage{amsthm}
\usepackage[pdftex,unicode]{hyperref}
\usepackage{amscd}
\usepackage[all,cmtip]{xy}
\usepackage{mathtools}
\usepackage{comment}
\usepackage{tikz-cd}

\def\Sl{\mathbb S}  \def\R{\mathbb R} \def\Q{\mathbb Q} 

\def\uu{\Upsilon}

\def\om{\omega}

\def\aut#1{\operatorname{Aut}(#1)}

\def\autp#1#2{\operatorname{Aut_{#1}}(#2)}
\def\be{\beta}

\def\al{\alpha}
\def\ga{\gamma}
\def\la{\lambda}
\def\ph{\varphi}
\def\de{\delta}

\def\ggg{\Gamma}
\def\tggg{\wh\ggg}
\def\gbm{\ggg^{\scriptscriptstyle{BM}}}
\def\god{\ggg^{\scriptscriptstyle{OD}}}

\def\nfs/{NSF}
\def\cdp/{CDP}
\def\cdpz/{CDP${}_0$}
\def\cnL#1{$(L_{#1})$}
\def\cnLl#1{$(L'_{#1})$}

\def\lleq{\preccurlyeq}
\def\ggeq{\succcurlyeq}

\def\cQ{\mathcal{Q}}

\def\cW{\mathcal{W}}
\def\cN{\mathcal{N}}
\def\cL{\mathcal{L}}
\def\cp{\mathcal{p}}
\def\cp{\mathcal{S}}
\def\s{{\mathfrak{Q}}}
\def\cG{\mathcal{G}}

\def\gr{{\mathfrak{R}}}

\def\uu{\Upsilon}

\def\ww{{\mathfrak{W}}}
\def\vv{{\mathfrak{V}}}

\def\vvbm{\vv_{BM}}

\def\pn{{\mathcal{N}}}

\def\sq#1#2{(#1)_{#2}}
\def\sqn#1{\sq{#1}{n\in\om}}
\def\sqnn#1{\sqn{#1_n}}

\def\markS#1{$(Q_{#1})$}

\def\llrrd#1#2{\sideset{^{{\rlap{$\scriptstyle{#1}$}\,\,}}}{{#2}}{\operatorname{\D}}}

\def\lrd#1#2{\llrrd{#1}{^{#2}}}
\def\rd#1{\D^{#1}}
\def\ld#1{\llrrd{#1}{}}

\def\lrdd#1#2#3{\llrrd{#1}{^{#2}_{#3}}}
\def\rdd#1#2{\D^{#1}_{#2}}
\def\ldd#1#2{\llrrd{#1}{_{#2}}}
\def\dd#1{\D_{#1}}

\def\lrds#1#2#3{\lrd{#1}{#2}(#3)}

\def\ds#1{\D(#1)}
\def\dps#1#2{\ds{#1;#2}}
\def\dpss#1{\dps{#1}{\s}}
\def\dpssp{\dpss{\cP}}

\def\cl#1{\overline{#1}}
\def\clx#1#2{\overline{#2}^{#1}}
\def\clv#1{\clx{v}{#1}}
\def\clh#1{\clx{h}{#1}}
\def\cld#1{\clx{d}{#1}}
\def\cla#1{\clx{a}{#1}}

\def\Int{\operatorname{Int}}

\def\bt{\operatorname{\beta}}

\def\ltu{\operatorname{Ls}}

\def\oms{\om^*}

\def\limp#1{{\textstyle \lim_{#1}}}

\def\es{\varnothing}
\def\Tau{{\mathcal T}}

\def\nom{{n\in\om}}

\def\sset#1{\{#1\}}
\def\fset#1{\{#1\}}
\def\set#1{\bbset#1\eeset}
\def\bbset#1:#2\eeset{\{#1\,:\,#2\}}

\def\bbsett#1:#2\eesett{\{#1\,:\,\text{#2}\}}

\def\iset#1{\ibbset#1\ieeset}
\def\ibbset#1:#2\ieeset{(#1)_{#2}}

\def\sod#1{\mathfrak{S}(#1)} 
\def\sse#1{\Expne {#1 \times #1}}

\def\tp{\Tau}
\def\tps{{\tp^*}}

\def\cP{{\mathcal P}}
\def\cB{{\mathcal B}}

\def\cB{{\mathcal B}}

\def\cU{{\mathcal U}}

\def\cR{{\mathcal R}}

\def\Exp#1{\operatorname{Exp}(#1)}
\def\Expne#1{\operatorname{{Exp}}_*(#1)}

\def\B{\mathbb D}
\def\til{\tilde}
\def\tD{{\til{\D}}}

\def\wt#1{\widetilde{#1}}

\def\eqdef{\coloneqq}

\def\cX{{\mathcal X}}
\def\cV{{\mathcal V}}

\def\nov#1{\operatorname{Nov}(#1)}
\def\dev#1{\operatorname{dev}(#1)}
\def\St{\operatorname{St}}
\def\st{\operatorname{st}}

\newcommand\restrA[2]{{
  \left.\kern-\nulldelimiterspace 
  #1 
  \vphantom{\big|} 
  \right|_{#2} 
  }}

\newcommand\restrB[2]{\ensuremath{\left.#1\right|_{#2}}}

\def\restr#1#2{\restrB{#1}{#2}}
\def\sorestr#1#2{\restr{#1}{#2}}

\def\pwr#1_#2{#1^{[#2]}}

\def\ddpsV#1#2{(#1)_{#2}}

\def\D{\Delta}
\def\term#1{{\it #1}}

\def\alp{{\al\in P}}

\def\ppl{{\mathsf{p}}}
\def\pp{{\mathsf{P}}}

\def\dddm#1(#2){N_{#1}(#2)}
\def\dddb#1(#2){B_{#1}(#2)}

\def\et(#1){}

\def\wh#1{\widehat{#1}}
\def\tg#1#2{\wh{#1#2}}
\def\et(#1){ (#1)}
\def\tbm #1{T_{BM}(#1)}
\def\tod #1#2{T_{OD}(#1,#2)}
\def\gbmt{\ggg^{\wh{\scriptscriptstyle{BM}}}}
\def\godt{\ggg^{\wh{\scriptscriptstyle{OD}}}}
\def\bitem#1,#2.{ $#2\nrightarrow #1$:\ }
\long\def\edemo#1\endedemo{{\rm #1}}

\newtheorem{statement}{Statement}[section]
\newtheorem{proposition}{Proposition}[section]
\newtheorem{theorem}{Theorem}[section]
\newtheorem*{theorem*}{Theorem}

\newtheorem*{lemma*}{Lemma}
\newtheorem{cor}{Corollary}[section]
\newtheorem{example}{Example}[section]

\theoremstyle{definition}
\newtheorem{definition}{Definition}[section]

\newtheorem{problem}{Problem}[section]
\theoremstyle{remark}
\newtheorem{note}{Remark}[section]
\newtheorem*{note*}{Remark}

\def\rarr{\Rightarrow}
\def\larr{\Leftarrow}
\def\lrarr{\Leftrightarrow}
\def\pppc{\pp_+^v}

\begin{document}
\begin{frontmatter}

\title{Classes of Baire spaces defined by semi-neighborhoods of the diagonal}

\author{Evgenii Reznichenko} 
\ead{erezn@inbox.ru}

\address{Department of General Topology and Geometry, Mechanics and  Mathematics Faculty, 
M.~V.~Lomonosov Moscow State University, Leninskie Gory 1, Moscow, 199991 Russia}

\begin{abstract}
With the help of semi-neighborhoods of the diagonal, classes of Baire spaces are defined: $\D$, $\D_h$ and $\D_s$ Baire spaces. These classes of spaces are studied with the help of topological games. They are useful in studying continuity in groups: paratopological $\D$-Baire groups, quasi-topological $\D_h$-Baire groups, and semitopological $\D_s$-Baire groups are topological groups.
\end{abstract}

\end{frontmatter}

\section{Introduction}\label{sec-intro}

A space $X$ is called Baire (nonmeager) if for any family $(U_n)_n$ of open dense subsets of $X$ the intersection $G=\bigcap_n U_n$ is dense in $X$ (not empty).

Baire spaces play an important role in mathematics. The strongest and most interesting results are obtained when additional constraints on the space are imposed, such as metrizability.

The concepts of a Baire space and nonmeager space  are closely related to each other.
\begin{theorem}\label{tintro1}
\begin{itemize}
\item[\rm(1)]
Baire spaces are nonmeager.
\item[\rm(2)]
An open subset of a Baire space is a Baire space.
\item[\rm(3)]
A space $X$ is  Baire if and only if every open subspace of $X$ is a nonmeager space.
\item[\rm(4)]
A space $X$ is nonmeager if and only if there exists an open nonempty Baire subspace of $X$.
\item[\rm(5)]
A homogeneous nonmeager space is a Baire space.
\end{itemize}
\end{theorem}
In the article, 
strengthening of Baire and nonmeager properties are constructed in parallel, and the relationships (1)--(5) between Baire and nonmeager spaces are checked.
The most important classes from the introduced classes are
$\D$-Baire, $\D_h$-Baire, $\D_s$-Baire and
$\D$-nonmeager, $\D_h$-nonmeager, $\D_s$-nonmeager spaces.

The main application of $\D$-Baire spaces is the continuity of operations in groups with the topology \cite{rezn2008,rezn2022-2008,moors2017}. The article \cite{Reznichenko2024rtg} is devoted to the study of the continuity of operations in 
topologized groups
that belong to the classes of spaces studied in this article.
%
The article \cite{Reznichenko2024rtg} is a revised part of the preprint \cite{rezn2022-1}.
There are two main research methods in this area. The first method is based on the exploitation of properties of compactness type. This approach is found in \cite{Montgomery1936,Reznichenko1994,Ravsky2001,Ravsky2002,RavskyReznichenko2002,arh-rezn2005,Ravsky2019,BanakhRavsky2020}.
The second method is based on a recursive procedure, which is best defined using modifications of the topological Banach--Mazur game.
This approach was used in \cite{Ellis1957,Zelazkoo1965,Brand1982,Pfister1985,Bouziad1993,Bouziad1996,Bouziad1996-2,KenderovKortezovMoors2001,arh2010,Reznichenko2024rtg}.
As noted in \cite{rezn2008,rezn2022-2008,moors2017}, spaces that are defined using topological games turn out to be $\D$-Baire in this approach. The largest subclass of $\D$-Baire spaces, which is defined using games, is described in \cite{rezn2008,rezn2022-2008} as the spaces $X$ for which in the game $G(OD,b)$ the player $\be $ has no winning strategy; in \cite{moors2017} such spaces are called Reznichenko spaces; in this article and \cite{rezn2022gbtg,Reznichenko2024rtg} these spaces are called $\god_{o,l}$-Baire.

Most theorems of the second type can be proved using the following steps: (1) using theorems on points of continuity of separately continuous functions of the Namioka type \cite{namioka1974}, properties of a group with a topology of the multiplication quasi-continuity type are proved; (2) properties of the topology of a group of type $\D$-Baire are proved; (3) using (1) and (2) we prove the continuity of operations in the group. Steps (1) and (2) are purely topological. Step (3) is algebro-topological. This article is about step (2), article \cite{Reznichenko2024rtg} is about step (3).

Let $X$ be a space, $P\subset X\times X$.
We call $P$ a \term{semi-neighbourhood of diagonal} if $x \in \Int P(x)$ for all $x\in X$ where $P(x)=\set{y\in X: (x,y)\in P}$. This article discusses several subclasses of Baire and nonmeager spaces, which are defined using semi-neighbourhoods of diagonal. For $\wt\D\in\sset{\D,\D_h,\D_s}$, a space $X$ is $\wt\D$-nonmeager ($\wt \D$-Baire) if for any semi-neighbourhood of the diagonal $P$ (and any nonempty open $O\subset X$)
there exists an open 
%
$W\subset X$ 
($W\subset O$), so that the condition ($\wt\D$) is satisfied:
\begin{itemize}
\item[($\D$)] $W\times W\subset \cl{P\cap(W\times W)}$;
\item[($\D_h$)] $W\subset \cl{\set{x: (x,y)\in P}}$ for all $y\in W$;
\item[($\D_s$)] $W\subset \cl{\set{x: \sset x \times W \subset P}}$.
\end{itemize}
The conditions ($\D$), ($\D_h$) and ($\D_s$) show how dense the set $S=P\cap(W\times W)$ is in $W\times W$. These conditions can be rewritten in the following form:
\begin{itemize}
\item[($\D$)] $S$ is dense in $W\times W$;
\item[($\D_h$)] $(W\times \sset x)\cap S$ is dense in $W\times \sset x$ for all $x\in W$;
\item[($\D_s$)] $M\times W\subset S$ for some set $M\subset W$ dense in $W$.
\end{itemize}
Clearly, $(\D_s)\rarr (\D_h)\rarr (\D)$.

The spaces that are here called ``$\D$-nonmeager'' were called ``$\D$-Baire'' in papers \cite{rezn2008,rezn2022-2008,tk2014,moors2017}.
I decided to rename the class of $\D$-Baire spaces because, as Proposition \ref{pdp2} shows, the $\D$-Baire property is a strengthening of the Baire property and the property of being $\D$-nonmeager is a strengthening of the property of being nonmeager. This does not cause any difficulties for applications of this property, since the $\D$-Baire and $\D$-nonmeager properties are the same for homogeneous spaces.

Some $\D$-Baire spaces can be found without games.
It is easy to check that a Baire metrizable space is $\D_s$-Baire. Let $M_n$ be the set of points $x\in X$ such that the ball
of radius $\frac 1{2^n}$ centered at $x$ is 
contained in
$P(x)$. Then $X=\bigcup_n M_n$ and the Baire property implies that $O\cap \Int \cl{M_n}$ is not empty for some $n$. Then if $W\subset O\cap \Int \cl{M_n}$ is an open set whose diameter is less than $\frac 1{2^n}$, then $W$ satisfies the condition $(\D_s)$.


In \cite{cdp2010} studies spaces, which in this article are called \cdp/-spaces (see Section \ref{sec-cdp}). It turns out these spaces are $\D_s$-Baire (Proposition \ref{pdsp1}). But it is unclear whether these spaces are $\god_{o,l}$-Baire (Problem \ref{pqe7}). It is also unclear whether there is an example \cdp/-space that does not contain a dense Baire metrizable space without additional set-theoretic assumptions (Problem \ref{pqe8}). Under assumption MA+$\lnot$CH the space $\R^{\om_1}$ is \cdp/-space (Example \ref{eqe6}).

This article is the second in a series of three articles, each of which is of independent interest.
In the first article \cite{rezn2022gbtg} we study topological games that arise in the study of groups with topology and the classes of spaces that are defined by these games. This article, the second article in the series, explores the properties of spaces that are defined using semi-neighbourhoods of diagonal, such as $\D$, $\D_h$ and $\D_s$ Baire and others. In the third paper, \cite{Reznichenko2024rtg}, we study groups that belong to the classes of spaces considered here.

In addition to the $\D$, $\D_h$ and $\D_s$ Baire spaces mentioned here, the article also introduces $\rd \ga$-Baire and $\lrd c\ga$-Baire spaces. Each $\D$-Baire space is $\rd \ga$-Baire and $\lrd c\ga$-Baire space. For $1<\ga\leq\om$ it is not known if there are $\rd \ga$-Baire and $\lrd c\ga$-Baire spaces that are not $\D$-Baire (Problem \ref{pqe9}).


The paper is organized as follows.

This Section \ref{sec-intro} is an introduction.
Section \ref{sec-defs} contains definitions and notation.

Sections \ref{sec-seminghbr} and \ref{sec-seminghbr2} introduces a general theory of topological properties, which are defined using semi-neighborhoods of the diagonal and normal square functors (\nfs/).
In Section \ref{sec-seminghbr}
$\D$, $\D_h$, $\D_s$ Baire and other similar classes of spaces are defined and their properties are studied.
Theorem 
\ref{t:nfs:1}
shows that the introduced classes of spaces are a correct strengthening of Baire and nonmeager properties, that is, conditions (1)--(5) of Theorem \ref{tintro1} are satisfied for these classes.
Section \ref{sec-seminghbr2} introduces even more subclasses of the class Baire spaces, which contain the class $\D$-Baire spaces.


Using this theory, \cite{Reznichenko2024rtg} defines and studies the algebraic-topological properties of groups with topology.
In addition, this theory is applicable to other \nfs/ that are useful in the study of groups with topology.

Section \ref{sec-bmg} deals with the topological games from \cite{rezn2022gbtg}, which are used to explore the spaces in the previous section.

Section \ref{sec-gd} establishes a relationship between the space classes from the previous two sections. The theorem \ref{tgd1} specifies specific classes of spaces with properties $\D$, $\D_h$, $\D_s$ Baire.

Section \ref{sec-cdp} introduces and studies \cdp/-spaces, which were studied in \cite{cdp2010} and defined by \cdp/-Baire spaces. It is proved that \cdp/-Baire spaces are $\D_s$-Baire (Proposition \ref{pdsp1}).

Section \ref{sec-tg} proves Theorem \ref{ttg1} about topologized groups, which is used to distinguish the introduced classes of spaces.

Section \ref{sec-qe} builds examples that distinguish the classes of spaces being studied. Problems are set.

\section{Definitions and notation}\label{sec-defs}

The sign $\eqdef$ will be used for equality by definition.


The family of all subsets of a set $X$ is denoted by $\Exp X$.
The family of all nonempty subsets of a set $X$ is denoted by $\Expne X$: $\Expne X\eqdef\Exp X \setminus \fset{\es}$.

If $B$ is a subset of a set $A$, then we denote by $B^c=A\setminus B$ the complement to $A$. We use this notation in situations where it is clear from the context which set $A$ is meant.

An \term{indexed set} $x=\iset{x_\al:\al\in A}$ is a function on $A$ such that $x(\al)=x_\al$ for $\al\in A$.
If the elements of an indexed set $\cX=\iset{X_\al: \al\in A}$ are themselves sets, then $\cX$ is also called an \term{indexed family of sets}; such an $\cX$ is a function on $A $ such that $\cX(\al)=X_\al$ for $\alp$.


We denote by $\aut X$ the set of all homeomorphisms of the space $X$ onto itself.


Let $M\subset X$. If $M$ is the union of a countable number of nowhere dense sets, then $M$ is called a \term{meager} set. Nonmeager sets are called \term{of the second Baire category} sets.
A subset of $M$ is called \term{residual} (comeager) if $X\setminus M$ is a meager set.

A topological space is called \term{meagre} (respectively, \term{nonmeagre}) if it is a meagre (respectively, nonmeagre) subset of itself. 

A family $\nu$ of nonempty subsets of $X$ is called a \term{$\pi$-net} if for any open nonempty $U\subset X$ there exists an $M\in\nu$ such that $M\subset U $.

A $\pi$-network consisting of open sets is called a \term{$\pi$-base}.

A subset $U\subset X$ is called \term{regular open} if $U=\Int{\cl U}$.

A space $X$ is called \term{quasi-regular} if for every nonempty open $U\subset X$ there exists a nonempty open $V\subset X$such that $\cl V\subset U$.

A space $X$ is called \term{semiregular} if $X$ has a base consisting of regular open sets.

A space $X$ is called \term{$\pi$-semiregular} \cite{Ravsky2001} (or nearly regular \cite{Ameen2021}) if $X $ is a $\pi$-base consisting of regular open sets.

For a cardinal $\tau$, a set $G\subset X$ is called a \term{set of type $G_\tau$} if $G$ is the intersection of $\tau$ open sets. A space $X$ is called \term{absolute $G_\tau$} if $X$ is of type $G_\tau$ in some compact extension.

A space \term{$X$ is regular at a point $x\in X$} if for any neighborhood $U$ of the point $x$ there exists a neighborhood $V\ni x$ such that $\cl V\subset U$.

A space \term{$X$ is semiregular at a point $x\in X$} if there is a base at the point $x$ consisting of regular open sets.

A space $X$ is \term{feebly compact} if any locally finite family of open sets is finite.

For $\ga\subset \Exp X$ and $x\in X$ we denote
\begin{align*}
\St(x,\ga)&\eqdef\set{U\in\ga: x\in \ga},
&
\st(x,\ga)&\eqdef\bigcup\St(x,\ga).
\end{align*}

A space $X$ is called \term{developable} if there exists a sequence of open covers $(\ga_n)_{n\in\om}$ such that for any $x\in X$ the family $\st(x,\ga_n)$ is a base at the point $x$.

A family $\cB$ of open nonempty sets in $X$ is called an \term{outer base of $M\subset X$} if $M\subset U$ for each $U\in\cB$ and for each open $W\supset M$ there exists a $U\in \cB$ such that $M\subset U\subset W$.

If $\sqnn M$ is a sequence of subsets of a space $X$, then the set
\begin{align*}
\ltu_{n\in\om} M_n \eqdef \{x\in X:\ &|\set{n\in\om: U\cap M_n\neq\es}|=\om
\\
&\text{ for any neighborhood }U\text{ of the point }x\}
\end{align*}
is called the \term{upper limit of the sequence of sets $\sqnn M$} \cite[Section 29]{KuratowskiBook1}.

If $\sqnn x$ is a sequence of points in a space $X$, then we denote
\begin{align*}
\ltu_{n\in\om} x_n & \eqdef \ltu_{n\in\om} \sset{x_n}.
\end{align*}

We denote by $\bt \om$ the space of ultrafilters on $\om$, the Stone-\v{C}ech extension of the discrete space $\om$. We denote by $\oms=\bt\om\setminus \om$ the set of nonprincipal ultrafilters.

Let $\sqnn x$ be a sequence of points in a space $X$ and let $p\in \oms$ be a nonprincipal ultrafilter. A point $x\in X$ is called the \term{$p$-limit} of a sequence $\sqnn x$ if $\set{n\in\om: x_n\in U}\in p$ for any neighborhood $U$ of $x$. We will write $x=\limp p x_n=\limp p \sqnn x$ for the $p$-limit $x$.

A \term{topologized group} is a triple $(G, \cdot, \tau )$, where $(G,\cdot)$ is a group and $(G, \tau )$ is a topological space \cite{cdp2010}. The neutral element of $(G, \cdot)$ is denoted by $e$.

\section{Space classes defined by semi-neighborhoods of the diagonal}\label{sec-seminghbr}

Let $X$ be a set, $P,Q\subset X\times X$, and $\cP\subset \Exp{X\times X}$. We denote
\begin{align*}
P(x)&\eqdef\set{y\in X: (x,y)\in P},
\\
P(M)&\eqdef\set{y\in X: \text{ there exists an }x\in M \text{ such that }(x,y)\in P}=\bigcup_{x\in M}P (x),
\\
\sorestr PM &\eqdef P\cap (M\times M),
\\
\sorestr \cP M &\eqdef \set{\sorestr PM: P\in \cP},
\\
P^{-1}&\eqdef\set{(x,y):(y,x)\in P},
\\
P\circ Q&\eqdef\set{(x,y):\text{there is a }z\in X\text{ such that }(x,z)\in P\text{ and }(z,y)\in Q}
\end{align*}
for $x\in X$ and $M\subset X$.

\begin{definition}
Let $X$ be a space, $P\subset X\times X$. We call a set $P$ \term{semiopen} if every $P(x)$ is open.
We call $P$ a \term{semi-neighbourhood of diagonal} if $x \in \Int P(x)$ for all $x\in X$.
We denote by $\sod X$ the family of all semi-neighborhoods of the diagonal of the space $X$.
\end{definition}

\begin{note}
In \cite{rezn2008}, in the definition of a semi-neighbourhood of a diagonal, it was assumed that the set $P(x)$ is open, that is, $P$ was assumed to be semiopen. The concept of a semiopen semi-neighborhood of a diagonal essentially coincides with the concept of neighbourhood assignment \cite{mtw2007}: a mapping $P:X\to \tp$ is called a \term{neighbourhood assignment} if $x\in P(x)$ and $ P(x)$ is an open neighborhood of $x\in X$.

A set $P$ is a semi-neighborhood of the diagonal if and only if $P$ contains some semiopen semi-neighborhood of the diagonal.
\end{note}


A functor $\s$ that associates a topological space $X$ with a family $\s(X)\subset \sse X$ of nonempty subsets of $X\times X$ is called a \term{square functor}. A square functor $\s$ is called a \term{normal square functor} (\nfs/) if the following conditions are satisfied:
\begin{itemize}
\item[\markS1] if $f:X\to Y$ is a homeomorphism, then
\[
\s(Y)=\set{(f\times f)(P): P\in \s(X)};
\]
\item[\markS2] if $U$ is an open nonempty subset of $X$, then $\s(U)=\sorestr {\s(X)}U$;
\item[\markS3] if $S\in \s(X)$, $S\subset Q \subset X\times X$, then $Q\in \s(X)$;
\item[\markS4] if $S\in \s(X)$, then $\cl S = X\times X$.
\end{itemize}

\begin{definition}\label{d:nfs:1}
Let $\s$ be a normal square functor, $X$ be a space,
and  $\cP\subset \sod X$. We call $X$
\begin{itemize}
\item a \term{$\dpss \cP$-nonmeager}  space if for any $P\in\cP$ there exists an open nonempty $V\subset X$ such that $\sorestr PV\in \s(V)$;
\item \term{$\dpss \cP$-Baire} space if for any $P\in\cP$ and any open nonempty $U\subset X$ there exists an open nonempty $V\subset U$ such that $\sorestr PV\in \s(V)$.
\end{itemize}
We call $X$
\begin{itemize}
\item a \term{$\ds\s$-nonmeager}  space if $X$ is $\dpss{\sod X}$-nonmeager;
\item a \term{$\ds\s$-Baire}  space if $X$ is $\dpss{\sod X}$-Baire.
\end{itemize}
\end{definition}

In this section we study $\ds\s$-nonmeager and $\ds\s$-Baire spaces. Section \ref{sec-seminghbr2} studies $\dpss \cP$-nonmeager and $\dpss \cP$-Baire spaces for proper families $\cP\subset \sod X$. $\dpss \cP$-Baire spaces are also used to study topologized groups \cite{Reznichenko2024rtg}.

The following theorem is an analogue of Theorem \ref{tintro1}.

\begin{theorem}\label{t:nfs:1}
Let $\s$ be a normal square functor.
\begin{itemize}
\item[\rm(1)]
$\ds\s$-Baire spaces are $\ds\s$-nonmeager.
\item[\rm(2)]
An open subset of a $\ds\s$-Baire space is a $\ds\s$-Baire space.
\item[\rm(3)]
A space $X$ is  $\ds\s$-Baire if and only if every open subspace of $X$ is a $\ds\s$-nonmeager space.
\item[\rm(4)]
A space $X$ is $\ds\s$-nonmeager if and only if there exists an open nonempty $\ds\s$-Baire subspace of $X$.
\item[\rm(5)]
A homogeneous $\ds\s$-nonmeager space is a $\ds\s$-Baire space.
\end{itemize}
\end{theorem}
\begin{proof}
As with Theorem \ref{tintro1}, conditions (1), (2) and (3) are obvious.
Let $(X,\tp)$ be a space, and $\tps=\tp\setminus\sset\es$.
Let us denote
\begin{align*}
\cU &= \set{U\in\tps : U \text{ is $\ds\s$-Baire}},
&
W &= \bigcup \cU,
\\
\cV &= \set{V\in\tps : V \text{ is $\ds\s$-nonmeager}},
&
\cQ &= \tps\setminus \cV.
\end{align*}
If $\cU\neq\es$, then $W\in\cU$ and $\cU=\set{U\in\tps: U\subset W}$. The space $X$ is $\ds\s$-Baire if and only if $\cV$ is a $\pi$-base of $X$.

(4) This condition is equivalent to the condition: $X$ is $\ds\s$-nonmeager $\lrarr$ $\cU\neq\es$.
($\larr$) Obviously.
($\rarr$) Let us assume the opposite, that is, $\cU=\es$. From (3) it follows that $\cQ$ is a $\pi$-base of $X$.
Let $\ga$ be a maximal disjoint subfamily of the family $\cQ$. Then $\bigcup\ga$ is dense in $X$ and $F=X\setminus \bigcup \ga$ is closed and nowhere dense in $X$. For each $G\in\ga$ we fix $P_G\in\sod G$,  which ensures that $G$ is not a $\ds\s$-nonmeager space, that is, $\sorestr {P_G}U\notin \s(U)$ for every open nonempty $U\subset G$. Let us put
\[
P = (F\times X) \cup \bigcup_{G\in\ga} P_g.
\]
Then $P\in \sod X$ and $\sorestr {P}U\notin \s(U)$ for each $U\in\tps$. Therefore $X$ is not $\ds\s$-nonmeager. Contradiction.

(5) From (4) it follows that $\cU\neq\es$. Note that $f(U)\in\cU$ for each $U\in\cU$ and $f\in\aut X$. Since $X$ is a homogeneous space, then $W=\bigcup \cU = X$. Therefore, $X\in\cU$, that is, $X$ is $\ds\s$-Baire.
\end{proof}

For $S\subset X\times X$, we denote 
\begin{align*}
\clv S&=\bigcup_{x\in X}\cl{\sset x \times S(x)},
&
\clh S&=\bigcup_{y\in X}\cl{ S^{-1}(y) \times  \sset y}.
\end{align*}
Note that 
\begin{align*}
\clv S &= \left(\clh{S^{-1}}\right)^{-1},
&
\clh S &= \left(\clv{S^{-1}}\right)^{-1}.
\end{align*}
It is directly verified that
\begin{align}\label{eq:nfs:1}
\clv S(x)&=\cl{S(x)},
&
\clh S(x) &= \bigcap_{U\in \tp_x} S(U),
&
\cl S(x) &= \bigcap_{U\in \tp_x} \cl{S(U)},
\end{align}
where $x\in X$ and $\tp_x$ are the family of open neighborhoods of the point $x$.

Let us define the \nfs/s which we will use. Let $k\in\sset{d,h,v,s,a}$, $X$ be a space and $S\subset X\times X$. Then $S\in \s_k(X)$ if and only if the  corresponding  condition  \cnL k holds:
\begin{itemize}
\item[\cnL d] $\cl S=X\times X$;
\item[\cnL v] $\clv S=X\times X$;
\item[\cnL h] $\clh S=X\times X$;
\item[\cnL s] $M\times X\subset S$ for some $M\subset X=\cl M$;
\item[\cnL a] $S=X\times X$.
\end{itemize}
In the diagram below we indicate the relationship between the entered conditions.
\[
\begin{tikzcd}
(L_a) \ar[r]\ar[rrd] &(L_s) \ar[r] & (L_h) \ar[r]& (L_d)
\\
&& (L_v) \ar[ru] &
\end{tikzcd}
\]

\begin{definition}\label{d:dp1:1}
Let  $X$ be a space. 
For ${k}\in\sset{d,v,h,s,a}$, we call $X$ 
\begin{itemize}
\item a \term{$\D_{{k}}$-nonmeager}  space if $X$ is $\ds{\s_{{k}}}$-nonmeager;
\item a \term{$\D_{{k}}$-Baire}  space if $X$ is $\ds{\s_{{k}}}$-Baire.
\end{itemize}
We call $X$ 
\begin{itemize}
\item a \term{$\D$-nonmeager}  space if $X$ is $\D_d$-nonmeager;
\item a \term{$\D$-Baire}  space if $X$ is $\D_d$-Baire.
\end{itemize}
\end{definition}

The definition of $\D$-, $\D_h$-, and $\D_s$-nonmeager (Baire) spaces was also given in the Introduction (Section \ref{sec-intro}).

For $S\subset X\times X$, we denote $\cld S=\cl S$ and $\cla S=S$. Then we can reformulate the definition of $\D_k$-Baire and $\D_k$-nonmeager spaces for $k\in\sset{d,h,v,a}$:
a space $X$ is $\D_k$-nonmeager ($\D_k$-Baire) if for any $P\in\sod X$ (and open nonempty $V\subset X$) there exists open nonempty $U\subset X$ ($U\subset V$) such that $U\times U \subset \clx kP$.

\begin{proposition} \label{p:dp1:1}
Let $X$ be a regular space. 
\begin{itemize}
\item[{\rm ($1$)}] If $X$ is $\D_v$-nonmeager  then $X$ is $\D_a$-nonmeager.
\item[{\rm ($2$)}] If $X$ is $\D_v$-Baire then $X$ is $\D_a$-Baire.
\end{itemize}
\end{proposition}
\begin{proof}
Let $P\in \sod X$. Since $X$ is regular then there exists $Q\in \sod X$ such that $\cl{Q(x)}\subset P(x)$ for any $x\in X$.

(1)
Since $X$ is $\D_v$-nonmeager, then there is a nonempty open $U\subset X$ such that $U\subset \cl{Q(x)}$ for any $x\in U$. Then $U\times U\subset P$.

(2)
Since $X$ is $\D_a$-Baire, then for any nonempty open $V\subset X$ there is a nonempty open $U\subset V$ such that $U\subset \cl{Q(x)}$ for any $x\in U$. Then $U\times U\subset P$.
\end{proof}

In the diagram below we indicate the relationship between the entered classes.
An arrow from $A$ to $B$ means that any $A$-nonmeager space is $B$-nonmeager and any $A$-Baire space is $B$-Baire.
\[
\begin{tikzcd}
\D_a \ar[r]\ar[rrd] &D_s \ar[r] & \D_h \ar[r]& \D 
\\
&& 
\D_v \ar[ru]
\arrow[bend left=20]{llu}[near end]{\text{for regular spaces}}
&
\end{tikzcd}
\]

\begin{proposition} \label{p:dp1:2}
Any $\D_h$-nonmeager space is a nonmeager space.
\end{proposition}
\begin{proof}
Let us prove that if $X$ is a meager space, then $X$ is not $\D_h$-nonmeager.
There is an increasing sequence $\sqnn F$ of closed nowhere dense subsets of the space $X$, so that $F_0=\es$ and $\bigcup_\nom F_n=X$. Let us denote $U_n=X\setminus F_n$. Let $f(x)=\max \set{\nom: x \in U_n}$ and
\[
P=\set{(x,y)\in X\times X: y\in U_{f(x)}}.
\]
Then $P$ is semi-neighborhood of the diagonal and
$P^{-1}(y)=\set{x\in X: y\in U_{f(x)}=\set{x\in X: f(x) \leq f(y)} }=F_{f(y)+1}$.
Therefore, $\clh P=P$. 

Let $U\subset X$ is a nonempty open set. Let $m=\min {\nom: F_n\cap U\neq\es}$ and $x\in U_m\cap U$. Then $P(x)\subset U_m$. Since $U\not\subset U_m$, then $U\not\subset P(x)$ and, therefore, $U\times U\not\subset P=\clh P$. Therefore, $X$ is not $\D_h$-nonmeager.
\end{proof}

\begin{statement} \label{adp1}
Let $X$ be a quasi-regular space.
A space $X$ is meager if and only if there exists a sequence $\ga_n$ of families of nonempty open sets for which the following conditions are satisfied:
\begin{enumerate}
\item $\ga_0=\sset X$;
\item the family $\ga_n$ is disjoint, $\cl{\bigcup\ga_n}=X$ and $\ga_{n+1}$ refines $\ga_n$;
\item $\cl U\cap \cl V=\es$ for different $U,V\in\ga_n$;
\item if $\nom$, $V\in\ga_{n+1}$, $U\in\ga_{n}$ and $V\cap U\neq\es$, then $\cl V\subset U$;
\item the family $\ga=\bigcup\limits_{n\in \om}\ga_n$ is pointwise finite.
\end{enumerate}
\end{statement}
\begin{proof}
 If $\ga_n$ with the listed properties exist, then $F_n=X\setminus\bigcup\ga_n$ is nowhere dense, $X=\bigcup_{n\in \om}F_n$, and hence $X$ is a meager space.

Suppose that $X$  is a meager space. Then there exist nowhere dense closed sets $F_n\subset X$, $n\in\om$ such that $X=\bigcup_{n\in \om}F_n$. We can assume that $F_0=\es$ and $F_n\subset F_{n+1}$ for $n\in\om$. We construct $\ga_n$ by induction on $n$.

We set $\ga_0=\sset X$.
Suppose that $\ga_n$ has been built. Let $U\in \ga_n$. It follows from Proposition \ref{psonbd8} that there exists a disjoint family $\nu_U\subset\tps$ such that
\begin{enumerate}
\item $\bigcup \nu_U\subset U\setminus F_{n+1}\subset\cl{\bigcup \nu_U}$;
\item $\cl V\subset U\setminus F_{n+1}$ for $V\in\nu_U$;
\item $\cl W\cap \cl V=\es$ for different $W,V\in\nu_U$.
\end{enumerate}
Let $\ga_{n+1}=\bigcup\set{\nu_U:U\in\ga_n}$.
\end{proof}

\begin{statement} \label{adp2-1}
Let $X$ be a meager quasi-regular space.
Then there is $P\in \sod X$  such that $P$ is nowhere dense in $X\times X$.
\end{statement}
\begin{proof}
Take $\ga_n$, $n\in\om$, as in Statement \ref{adp1}. We set
\begin{align*}
f(x)&=\max \set{m\in\om: x\in \bigcup \ga_m}
\end{align*}
for $x\in X$.
We set
\begin{align*}
P=\{ (x,y)\in X\times X :\
& \text{there exists a }U\in \ga_l,\text{ where }l=f(x),
\\
& \text{such that }x,y\in U
\}.
\end{align*}
Let $U,V\in\tps$. Let $y\in V$. Let $m=f(y) + 1$. Take $W\in \ga_m$ such that $U'=W\cap U\neq\es$. Take $U''\in \ga_{m}$ such that $W\subset U''$. Then $P(U')\subset U''$ and $y\notin \cl{U''}$. Let $V'=V\setminus \cl{U''}$. Then $\es\neq U'\times V'\subset U\times V$ and $P\cap U'\times V'=\es$.
\end{proof}

The following statement follows from Statement \ref{adp2-1}.

\begin{proposition} \label{p:dp1:3}
Any quasi-regular $\D$-nonmeager space is a nonmeager space.
\end{proposition}

There is $\D_v$-Baire meager $T_1$ space (Example \ref{eqe10}).

The following statements follow from Propositions \ref{p:dp1:2}, \ref{p:dp1:3} and Theorems \ref{tintro1}, \ref{t:nfs:1}.

\begin{theorem}\label{t:dp1:1}
Let  $X$ be a space. 
\begin{itemize}
\item
A space $X$ is a nonmeager if $X$ belongs to any of the following classes of spaces:
$\D_h$-nonmeager; $\D_s$-nonmeager; $\D_a$-nonmeager; quasi-regular $\D_v$-nonmeager; quasi-regular $\D$-nonmeager.
\item[\rm(2)]
A space $X$ is a Baire if $X$ belongs to any of the following classes of spaces:
$\D_h$-Baire; $\D_s$-Baire; $\D_a$-Baire; quasi-regular $\D_v$-Baire; quasi-regular $\D$-Baire.
\end{itemize}
\end{theorem}

\section{More classes defined by semi-neighborhoods of the diagonal}\label{sec-seminghbr2}

In this section introduces subclasses of the class Baire spaces, which contain the class $\D$-Baire spaces.

\subsection{Semi-neighborhoods of the diagonal}\label{sec-sonbd}

The following statement follows from the equalities \ref{eq:nfs:1}.

\begin{proposition} \label{psonbd1}
Let $X$ be a space, $P\subset X\times X$, $R\in \sod X$. Then
\begin{align*}
\clh P &= \bigcap \set{{Q\circ P}: Q\in \sod X},
&
\clh P & \subset R\circ P,
\\
\cl P &= \bigcap \set{\clv{Q\circ P}: Q\in \sod X},
&
\cl P & \subset \clv{R\circ P}.
\end{align*}

\end{proposition}

On $\Exp{\sod X}$ we introduce an order relation and an equivalence relation related to it. For $\cP,\cQ\subset \sod X$ we put
\begin{itemize}
\item
$\cP\lleq \cQ$ if and only if for any $Q\in \cQ$ there exists a $P\in\cP$ such that $P\subset Q$;
\item
$\cP\sim \cQ$ if and only if $\cP\lleq \cQ$ and $\cQ\lleq \cP$.
\end{itemize}
Note that if $\cP\supset \cQ$, then $\cP\lleq \cQ$.
For $\cP\subset \sod X$ we denote
\[
\pp^e(\cP)\eqdef\set{R\in \sod X: P\subset R,\text{ for some }P\in\cP}.
\]

\begin{proposition} \label{psonbd2}
Let $X$ be a space.
For $\cP,\cQ\subset \sod X$
\begin{itemize}
\item
$\cP\lleq \cQ$ if and only if $\pp^e(\cP)\supset \pp^e(\cQ)$;
\item
$\cP\sim \cQ$ if and only if $\pp^e(\cP) = \pp^e(\cQ)$.
\end{itemize}
\end{proposition}

For $\cP\subset \sod X$ we denote
\begin{align*}
\pp(\cP)&\eqdef\cP,
&
\pp_+(\cP)&\eqdef\set{Q\circ P: Q\in\sod X,\ P\in \cP},
\\
\pp^v(\cP)&\eqdef\set{\clv P: P\in\cP},
&
\pp^c(\cP)&\eqdef\set{\cl P: P\in\cP},
\\
\pppc(\cP)&\eqdef \pp^v(\pp_+(\cP)).
\end{align*}

Proposition \ref{psonbd1} implies the following assertion.

\begin{proposition} \label{psonbd3}
Let $X$ be a space.
For $\cP \subset \sod X$
\begin{itemize}
\item
$\cP\lleq \pp_+(\cP)\lleq \pppc(\cP)$,
\item
$\cP\lleq \pp^v(\cP)\lleq \pp^c(\cP)\lleq \pppc(\cP)$.
\end{itemize}
\end{proposition}

Let $\la$ be an ordinal and let $P_\al\subset X\times X$ for $\al<\la$. By induction on $\be<\ga$ we define sets $\ppl_\be(\cp_\be),\ppl^c_\be(\cp_\be)\subset X\times X$, where $\cp_\be= (P_\al)_{\al<\be}$.

For $\be=0$, we set
\[
\ppl_0(\cp_0)\eqdef\ppl_0^c(\cp_0)\eqdef\es.
\]

For $\be=1$, we set
\[
\ppl_1(\cp_1)\eqdef P_0,
\qquad
\ppl_1^c(\cp_1)\eqdef \cl{P_0}.
\]

For $1<\be\leq \ga$, we set
\begin{align*}
\ppl_\be(\cp_\be)&\eqdef\begin{cases}
\bigcup_{\al<\be} \ppl_\al(\cp_\al)&\text{if $\be$ is a limit ordinal},
\\
\ppl_{\de}(\cp_{\de})\circ P_{\de}&\text{if }\be=\de+1,
\end{cases}
\\
\ppl^c_\be(\cp_\be)&\eqdef\begin{cases}
\cl{\bigcup_{\al<\be} \ppl^c_\al(\cp_\al)}&\text{if $\be$ is a limit ordinal},
\\
\cl{\ppl_{\de}^c(\cp_{\de})\circ P_{\de}}&\text{if }\be=\de+1.
\end{cases}
\end{align*}

For $\cP\subset \sod X$ we denote
\begin{align*}
\pp_\ga(\cP)&\eqdef\set{\ppl_\ga(\cp): \cp\in \cP^\ga},
&
\pp^c_\ga(\cP)&\eqdef\set{\ppl^c_\ga(\cp): \cp\in \cP^\ga}
\end{align*}
for $Y\subset X$.

Proposition \ref{psonbd3} and definitions imply the following assertion.

\begin{proposition} \label{psonbd4}
Let $X$ be a space, $\cP \subset \sod X$, and $1<\al<\be$ ordinals.
Then
\begin{align*}
\pp_+(\cP) \lleq \pp_2(\cP) \lleq &\pp_\al(\cP) \lleq \pp_\be(\cP) \lleq \pp_\be^c(\cP),
\\
\pp^c(\cP) = \pp_1^c(\cP) \lleq &\pp_\al^c(\cP) \lleq \pp_\be^c(\cP),
\\
\pp^c(\cP) \lleq \pppc(\cP) \lleq &\pp_2^c(\cP).
\end{align*}
\end{proposition}

\begin{proposition} \label{psonbd4+1}
Let $X$ be a space, $Y\subset X$, $\cP \subset \sod X$, and $\ga$ be an ordinal.
If $\til\pp\in\sset{\pp^v,\pp^c,\pp_+,\pppc,\pp_\ga,\pp_\ga^c}$, then
\[
\til\pp(\sorestr \cP Y)\lleq \sorestr {\til\pp(\cP)} Y.
\]
\end{proposition}

\begin{proposition} \label{psonbd5}
Let $X$ be a space, $\la$ be an ordinal, and $U\subset X$ be a nonempty open subset of $X$.
Then
\begin{align*}
\sorestr{\pp_\ga(\sod X)} U &\sim \pp_\ga(\sod X) ,
&
\sorestr{\pp^c_\ga(\sod X)} U
&\sim
\pp^c_\ga(\sod U).
\end{align*}
\end{proposition}

\begin{definition} Let $(X,\tp)$ be a space, $\tps=\tp\setminus \sset\es$ and $\cP\subset \sod X$.
We call the family $\cP$ \term{partible} if the  following condition holds:
\begin{itemize}
\item[$(P)$] Let $A$ be an index set, a family $\set{U_\al: \al\in A}\subset \tps$ be disjoint and $\set{P_\al: \al\in A}\subset\cP$. Then there are a $P\in\cP$ and a family of open sets $\set{V_\al: \al\in A}\subset \tps$ such that for every $\al\in A$ the following conditions are satisfied:
\begin{itemize}
\item[$(1)$] $V_\al\subset U_\al \subset \cl{V_\al}$;
\item[$(2)$] $P(V_\al)\subset U_\al$;
\item[$(3)$] if $x\in V_\al$, then $P(x)\subset P_\al(x)$.
\end{itemize}
\end{itemize}
\end{definition}

The definition  implies the following assertion.

\begin{proposition} \label{psonbd6}
For any space $X$ the family $\sod X$ is partible.
\end{proposition}

\begin{proposition} \label{psonbd7}
Let $X$ be a space, $\cP,\cQ\subset \sod X$ and $\cR=\set{P\circ Q: P\in\cP\text{ and }Q\in\cQ}$. If $\cP$ and $\cQ$ are partible families, then $\cR$ is a partible family.
\end{proposition}
\begin{proof}
Let $\tp$ be the topology of $X$, and let $\tps=\tp\setminus \sset\es$.

Let $A$ be an index set, let  $\set{U_\al: \al\in A}\subset \tps$ be a disjoint family, and let $\set{R_\al: \al\in A}\subset\cR$, where $R_\al=P_\al\circ Q_\al$ for $\al\in A$. Since $\cP$ is a partible family, there exists a $P\in\cP$ and a family of open sets $\set{W_\al: \al\in A}\subset \tps$ such that for every $\al\in  A$ the following conditions are met:
 $W_\al\subset U_\al \subset \cl{W_\al}$;
$P(W_\al)\subset U_\al$;
if $x\in W_\al$, then $P(x)\subset P_\al(x)$.
Since $\cQ$ is a partible family, there exists a $Q\in\cQ$ and a family of open sets $\set{V_\al: \al\in A}\subset \tps$ such that for every $\al\in  A$ the following conditions are met:
 $V_\al\subset W_\al \subset \cl{V_\al}$;
 $Q(V_\al)\subset W_\al$;
 if $x\in V_\al$, then $Q(x)\subset Q_\al(x)$.
 
 Let $R=P\circ Q$. Then for each $\al\in A$ the following conditions are satisfied:
 $V_\al\subset U_\al \subset \cl{V_\al}$;
 $R(V_\al)\subset U_\al$;
 if $x\in V_\al$, then $R(x)\subset R_\al(x)$.
\end{proof}

\begin{proposition} \label{psonbd8}
Let $(X,\tp)$ be a quasi-regular space. For every open nonempty $U\subset X$, there exists a disjoint family $\ga$ of open subsets of $X$such that
\begin{itemize}
\item[$(1)$] $\bigcup\ga\subset U \subset \cl{\bigcup\ga}$;
\item[$(2)$] $\cl V\subset U$ for each $V\in \ga$.
\end{itemize}
\end{proposition}
\begin{proof}
Let $\tp$ be the topology of the space $X$, $\tps=\tp\setminus \sset\es$ and
\[
\nu=\set{V\in\tps:\cl V\subset U}.
\]
For $\ga$ any maximal disjoint subfamily of $\nu$ can be taken.
\end{proof}


\begin{statement} \label{asonbd1}
Let $(X,\tp)$ be a quasi-regular space and let $\tps=\tp\setminus \sset\es$, $\cP \subset \sod X$ be a partible family.
Let $A$ be an index set, $\set{U_\al: \al\in A}\subset \tps$ be a disjoint family, and  $\set{P_\al: \al\in A}\subset\cP$. Then there are $P\in\cP$ and a family of open sets $\set{V_\al: \al\in A}\subset \tps$ such that for every $\al\in A$ the following conditions are satisfied:
\begin{itemize}
\item[$(1)$] $V_\al\subset U_\al \subset \cl{V_\al}$;
\item[$(2)$] $\cl P(V_\al)\subset U_\al$;
\item[$(3)$] if $x\in V_\al$, then $P(x)\subset P_\al(x)$.
\end{itemize}
\end{statement}
\begin{proof}
We use Proposition \ref{psonbd8} and for each $\al\in A$ choose a disjoint family $\ga_\al\subset\tps$ such that $\bigcup\ga_\al\subset U_\al \subset \cl{\bigcup \ga_\al}$ and
$\cl U\subset U_\al$ for each $U\in \ga_\al$. We denote $\al(U)=\al$ and $P_U=P_\al$ for $U\in \ga_\al$.
We put $\ga=\bigcup_{\al\in A}\ga_\al$.
Since $\cP$ is a partible family, there exists a $P\in\cP$ and a family of open sets $\set{V_U: U\in\ga}\subset \tps$ such that for every $U\in\ga $ the following conditions are met:
\begin{itemize}
\item $V_U\subset U \subset \cl{V_U}$;
\item $P(V_U)\subset U$;
\item if $x\in V_U$ then $P(x)\subset P_U(x)$.
\end{itemize}
Let $V_\al=\bigcup\set{V_U: U\in \ga_\al}$ for $\al\in A$. By construction, conditions (1) and (3) are  satisfied.
Let us check (2). Let $x\in V_\al$. Then $x\in V_U$ for some $U\in \ga_\al$. We obtain
\[
P(V_U)\subset \cl{P(V_U)} = \cl U \subset U_\al.
\]
Proposition \ref{psonbd1} implies that $\cl P(x) \subset \cl{P(V_U)} \subset U_\al$.
\end{proof}

Propositions \ref{psonbd6}, \ref{psonbd7} and Statement \ref{asonbd1} imply the following assertion.
\begin{proposition} \label{psonbd9}
Let $X$ be a space and let $\cP\subset \sod X$ be a partible family.
\begin{itemize}
\item
The families $\pp_+(\cP)$ and $\pp_n(\cP)$ for $n<\om$ are partible.
\item
If $X$ is a quasi-regular space, then the families $\pp^v(\cP)$, $\pp^c(\cP)$, $\pppc(\cP)$ and $\pp_n^c( \cP)$ for $n<\om$ are partible.
\end{itemize}
\end{proposition}

\begin{definition} Let $X$ be a space and $\cP\subset \sod X$.
We denote
\begin{align*}
\autp{\cP}X \eqdef \{ f\in \aut X:\
&(f\times f)(P)\in \cP\text{ and}\\
&(f^{-1}\times f^{-1})(P)\in \cP\\
&\text{for all }P\in\cP
\}.
\end{align*}
We call a space $X$ \term{$\cP$-homogeneous} if for any $x,y\in X$ there exists an $f\in \autp{\cP}X$ such that $f(x)=y$ .
\end{definition}

The following proposition is a direct consequence of the definitions.

\begin{proposition} \label{psonbd10}
Let $X$ be a space, $\cP\subset \sod X$, $\al$ be an ordinal and
\[
\cQ\in\sset{\pp^e(\cP),\pp^v(\cP),\pp^c(\cP),\pp_+(\cP),\pppc(\cP),\pp_\al(\cP),\pp_\al^c(\cP)}.
\]
Then
\begin{itemize}
\item
$\autp \cP X\subset \autp \cQ X$;
\item
if $X$ is $\cP$-homogeneous, then $X$ is $\cQ$-homogeneous.
\end{itemize}
\end{proposition}

\begin{proposition} \label{psonbd11}
Let $X$ be a space, $\cP = \sod X$, $\al$ be an ordinal,
\[
\cQ\in\sset{\cP,\pp^e(\cP),\pp^v(\cP),\pp^c(\cP),\pp_+(\cP),\pppc(\cP),\pp_\al(\cP),\pp_\al^c(\cP)}.
\] 
Then
\begin{itemize}
\item
$\aut X = \autp \cQ X$;
\item
if $X$ is homogeneous, then $X$ is $\cQ$-homogeneous.
\end{itemize}
\end{proposition}

\subsection{Normal square functors}\label{sec-nfs}

Definitions  imply the following assertion.
\begin{proposition} \label{pnfs1}
Let $k\in\sset{d,h,v,s,a}$, $X$ be a space and $S\in \sse X$. Then $S\in\s_k(X)$ if and only if for any $x\in X$ and any neighborhood $U\subset X$ of $x$ the corresponding condition \cnLl k is satisfied:
\begin{itemize}
\item[\cnLl d] $\cl{S(U)}=X$;
\item[\cnLl v] $\cl{S(x)}=X$;
\item[\cnLl h] $S(U)=X$;
\item[\cnLl s] $S(z)=X$ for some $z\in U$;
\item[\cnLl a] $S(x)=X$.
\end{itemize}
\end{proposition}

We introduce an order relation on normal square functors. For \nfs/s $\s$ and $\gr$ $\s\lleq\gr$ is true if and only if $\s(X)\supset \gr(X)$ for any space $X$.

\begin{proposition} \label{pnfs2}
For any \nfs/ $\s$
\begin{align*}
&\s_d\lleq\s\lleq\s_a,
&
&\s_h\lleq\s_s.
\end{align*}
\end{proposition}

We introduce an increment operation on \nfs/: given an \nfs/ $\s$, we define an \nfs/ $\s^+$. Let $X$ be a space and let $Q\in\sse X$. Then $Q\in\s^+(X)$ if and only if $P\circ S\subset Q$ for some $P\in\sod X$ and $S\in\s(X)$.

\begin{proposition} \label{pnfs3}
\begin{align*}
& \s_v \lleq \s^+_d,
&
&\s_a = \s^+_s = \s^+_h.
\end{align*}
\end{proposition}

\subsection{$\dps \cP\s$-Baire spaces}\label{sec-dps}

The definitions imply the following assertion.

\begin{proposition} \label{pdps1}
Let $X$ be a space, $\s$, $\gr$ be \nfs/ and $\cP,\cR\subset \sod X$.
Suppose that $\s\ggeq \gr$ and $\cP\lleq\cR$.
\begin{itemize}
\item
If $X$ is a $\dps \cP\s$-Baire space, then $X$ is a $\dps \cR\gr$-Baire space.
\item
If $X$ is a $\dps \cP\s$-nonmeager space, then $X$ is a $\dps \cR\gr$-nonmeager space.
\end{itemize}
\end{proposition}

For $P\subset X\times X$ we denote
\begin{align*}
\cW(P,\s)&\eqdef \{V\in\tps: \sorestr PV\in \s(V)\},
\\
W(P,\s)&\eqdef\bigcup \cW(P,\s).
\end{align*}
For $\cP\subset \sod X$ we denote
\begin{align*}
\cB(\cP,\s)&\eqdef \{V\in\tps: \text{$V$ is $\dps{\sorestr{\cP}V}{\s}$-Baire}\},
\\
B(\cP,\s)&\eqdef\bigcup \cB(\cP,\s).
\end{align*}

The following three assertions are verified directly.

\begin{statement} \label{adps0}
$B(\cP,\s)\subset \cl{W(P,\s)}$ for $P\in\cP$.
\end{statement}

\begin{statement} \label{adps1}
For any $P\in\cP$ a space $X$ is
\begin{itemize}
\item $\dpss \cP$-nonmeager space if and only if $W(P,\s)\neq\es$,
\item $\dpss \cP$-Baire if and only if $\cl{W(P,\s)}=X$.
\end{itemize}
\end{statement}

\begin{statement} \label{adps2}
$W(P,\s)\cap V = W(\sorestr PV,\s)$ for $V\in\tps$.
\end{statement}

Statement \ref{adps1} implies the following proposition.

\begin{proposition} \label{pdps2}
Let $X$ be a space, $\s$ be an \nfs/ and $\cP\subset \sod X$.
If $X$ is $\dpssp$-Baire, then $X$ is $\dpssp$-nonmeager.
\end{proposition}

Statements  \ref{adps1} and \ref{adps2} imply the following proposition.
\begin{proposition} \label{pdps3}
Let $X$ be a space, $\s$ be an \nfs/ and $\cP\subset \sod X$.
\begin{itemize}
\item[\rm(1)]
A space $X$ is $\dpss \cP$-Baire if and only if $U$ is a $\dpss{\sorestr{\cP}U}$-nonmeager space for any nonempty open $U\subset X$.
\item[\rm(2)]
If $U\subset X$ is a nonempty open space and $U$ is a $\dpss{\sorestr{\cP}U}$-nonmeager space, then $X$ is $\dpss \cP$-nonmeager.
\end{itemize}
\end{proposition}

\begin{statement} \label{adps3}
If $P,Q\in\sod X$, then $W(Q,\s)\subset W(P\circ Q,\s^+)$.
\end{statement}
\begin{proof}
It suffices to check $\cW(Q,\s)\subset \cW(P\circ Q,\s^+)$. Let $V\in \cW(Q,\s)$. Then $\sorestr QV\in \s(V)$ and hence $\sorestr PV \circ \sorestr QV\in\s^+(X)$. Since $\sorestr PV \circ \sorestr QV\subset \sorestr {(P\circ Q)}V$, it follows that $\sorestr {(P\circ Q)}V \in \s^+(V)$.
\end{proof}

Statements  \ref{adps1} and \ref{adps3} imply the following proposition.
\begin{proposition} \label{pdps4}
Let $X$ be a space, $\s$ be an \nfs/ and $\cP\subset \sod X$.
\begin{itemize}
\item[\rm (1)] If $X$ is a $\dpss \cP$-nonmeager space, then $X$ is $\dps{\pp_+(\cP)}{\s^+}$-nonmeager space.
\item[\rm (2)] If $X$ is a $\dpss \cP$-Baire space, then $X$ is a $\dps{\pp_+(\cP)}{\s^+}$-Baire space space.
\end{itemize}
\end{proposition}

\begin{statement} \label{adps4}
If $P\in\sod X$, then $W(P,\s_d)\subset W(\cl P,\s_a)$.
\end{statement}
\begin{proof}
It suffices to check $\cW(P,\s_d)\subset \cW(\cl P,\s_a)$. Let $V\in \cW(P,\s_d)$. Then $\sorestr PV\in \s_d(V)$ and hence $V\times V\subset\cl{\sorestr PV}\subset \cl P$.
\end{proof}

Statements  \ref{adps1} and \ref{adps4} imply the following proposition.
\begin{proposition} \label{pdps5}
Let $X$ be a space and let $\cP\subset \sod X$.
\begin{itemize}
\item[\rm (1)] If $X$ is a $\dps \cP{\s_d}$-nonmeager space, then $X$ is a $\dps{\pp^c(\cP)}{\s_a}$-nonmeager space.
\item[\rm (2)] If $X$ is a $\dps \cP{\s_d}$-Baire space, then $X$ is a $\dps{\pp^c(\cP)}{\s_a}$-Baire space.
\end{itemize}
\end{proposition}


Proposition \ref{pdps5} implies the following assertion.

\begin{proposition} \label{pdps6}
Let $X$ be a space, $\cP\subset \sod X$, $\pp^c(\cP)=\cP$ and $\s$ be an \nfs/. Then
\begin{itemize}
\item[\rm (1)] $X$ is a $\dps \cP{\s}$-nonmeager space if and only if $X$ is a $\dps{\cP}{\s_a}$-nonmeager space;
\item[\rm (2)] $X$ is a $\dps \cP{\s}$-Baire space if and only if $X$ is a $\dps{\cP}{\s_a}$-Baire space.
\end{itemize}
\end{proposition}

Propositions \ref{psonbd3}, \ref{psonbd4},
\ref{pnfs2},
\ref{pnfs3},
\ref{pdps1},
\ref{pdps4},
\ref{pdps5}
and
\ref{pdps6}
imply the following assertion.

\begin{proposition} \label{pdps7}
Let $X$ be a space and let $\cP \subset \sod X$.
In the diagrams below, the arrow
\[
\ddpsV Fk \to \ddpsV Gl
\]
means that $F,G:\Exp{\sod X}\to \Exp{\sod X}$ are mappings, $k,l\in\sset{d,v,h,s,a}$ and the following  conditions hold:
\begin{itemize}
\item[\rm (1)] If $X$ is a $\dps{F(\cP)}{\s_k}$-nonmeager space, then $X$ is $\dps{G(\cP)}{\s_l}$-nonmeager space and
\item[\rm (2)] if $X$ is a $\dps{F(\cP)}{\s_k}$-Baire space, then $X$ is $\dps{G(\cP)}{\s_l}$-Baire space.
\end{itemize}

{
\def\xp#1{\ddpsV{\pp}#1}
\def\xpu#1#2{\ddpsV{\pp^{#1}}#2}
\def\xpl#1#2{\ddpsV{\pp_{#1}}#2}
\def\x#1#2#3{\ddpsV{\pp_{#1}^{#2}}#3}
\def\xc#1#2{\ddpsV{\pp_{#1}^{c}}#2}
\def\rr{\arrow[r,leftrightarrow]}

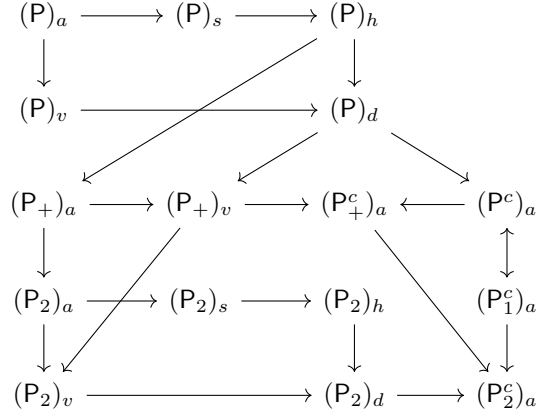
\[
\begin{tikzcd}
\xpu cv \rr&\xpu ca \rr&\xpu cs \rr&\xpu ch \rr&\xpu cd
\end{tikzcd}
\]

\[
\begin{tikzcd}
\xp a \arrow[r]\arrow[d] & \xp s \arrow[r] & \xp h \arrow[d] \arrow[ddll]& \\
\xp v \arrow[rr] && \xp d \arrow[dl]\arrow[dr]&\\
\xpl+a \arrow[r]\arrow[d] & \xpl+v \arrow[r] \arrow[ddl] & \x+va \arrow[ddr] &\xpu ca \arrow[l] \arrow [d,leftrightarrow]\\
\xpl 2a \arrow[r]\arrow[d] & \xpl 2s \arrow[r] & \xpl 2h \arrow[d] & \xc 1a \arrow[d]\\
\xpl 2 v \arrow[rr] && \xpl 2 d \arrow[r] & \xc 2a\\
\end{tikzcd}
\]

Let $1<\al<\be$ be ordinals.

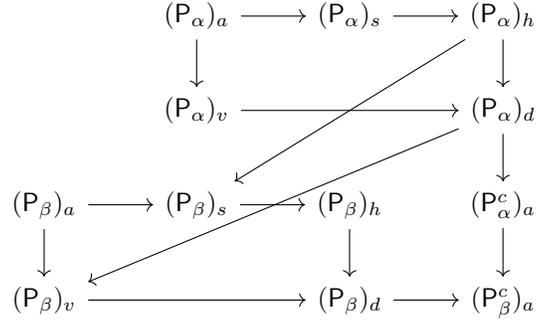
\def\xca#1{\xc{\al}{#1}}
\[
\begin{tikzcd}
\xca v \rr&\xca a \rr&\xca s \rr&\xca h \rr&\xca d
\end{tikzcd}
\]

\[
\begin{tikzcd}
&\xpl \al a \arrow[r]\arrow[d] & \xpl \al s \arrow[r] & \xpl \al h \arrow[ddll] \arrow[d] \\
&\xpl \al v \arrow[rr] && \xpl \al d \arrow[ddlll] \arrow[d]\\
\xpl \be a \arrow[r]\arrow[d] & \xpl \be s \arrow[r] & \xpl \be h \arrow[d] & \xc \al a \arrow[d]\\
\xpl \be v \arrow[rr] && \xpl \be d \arrow[r] & \xc \be a\\
\end{tikzcd}
\]
}

\end{proposition}

The following two assertions are easy to verify.

\begin{statement} \label{adps5}
Let $X$ be a space, $\s$ be an \nfs/ and $P\in\sod X$.
If $f\in\aut X$, then
\[
f(W(P,\s))=W((f\times f)(P),\s).
\]
\end{statement}

\begin{statement} \label{adps6}
Let $X$ be a space, $\s$ be an \nfs/ and $\cP\subset\sod X$.
If $f\in\autp \cP X$ then
\[
f(B(\cP,\s))=B(\cP,\s).
\]
\end{statement}

\begin{statement} \label{adps7}
Let $X$ be a space, $\cP\subset \sod X$ and $\s$ be an \nfs/. If $X$ is a $\cP$-homogeneous space and $B(\cP,\s)\neq\es$, then $X$ is $\dpss \cP$-Baire.
\end{statement}
\begin{proof}
Let $U=B(\cP,\s)$. Let us show that $U=X$. Take an arbitrary $y\in X$. Let $x\in U$. Since the space $X$ is $\cP$-homogeneous, there exists an $f\in\autp \cP X$ such that $f(x)=y$. From Statement \ref{adps6} it follows that $f(U)=U$. Hence $y\in U$. Proved $U=X$.

Let $P\in\cP$. It follows from Statement \ref{adps0} that $U\subset \cl{W(P,\s)}$. Hence $\cl{W(P,\s)}=X$. It follows from Statement \ref{adps1} that $X$ is $\dpss \cP$-Baire.
\end{proof}


\begin{statement} \label{adps8}
Let $X$ be a space, $\cP\subset \sod X$ be a partible family and $\s$ be an \nfs/.
The space $X$ is $\dpss \cP$-nonmeager if and only if $B(\cP,\s)\neq\es$.
\end{statement}
\begin{proof}
Let $\tp$ be the topology of $X$ and let $\tps=\tp\setminus \sset\es$.
Propositions \ref{pdps2} and \ref{pdps3} imply that if $U$ is a $\dpss{\sorestr{\cP}U}$-Baire space for some $U\in\tps$, then $X $ is $\dpss \cP$-nonmeager.

Let us prove the converse assertion. Assume that, on the contrary, $B(\cP,\s)=\es$. Then $U$ is not a $\dpss{\sorestr{\cP}U}$-Baire space for any $U\in\tps$. Proposition \ref{pdps3} implies that $V$ is not $\dpss{\sorestr{\cP}V}$-nonmeager for some $V\in\tps$, $V\subset U$. We set
\[
\nu=\set{V\in\tps: V\text{ is not $\dpss{\sorestr{\cP}V}$-nonmeager}}.
\]
The family $\nu$ is a $\pi$-base of $X$. Let $\ga\subset \nu$ be a maximal disjoint subfamily of the family $\nu$. Then $\cl{\bigcup\ga}=X$. It follows from Statement \ref{adps1} that for each $V\in\ga$ there exists a $P_V\in\cP$ for which
\[
W(P_V,\s)\cap V= W(\sorestr{P_V}V,\s)=\es.
\]
Since $\cP$ is a partible family, there exist $P\in\cP$ and $W_V\in\tps$ for $V\in\ga$ such that the   following conditions hold:
\begin{itemize}
\item $W_V\subset V \subset \cl{W_V}$;
\item $P(W_V)\subset V$;
\item if $x\in W_V$ then $P(x)\subset P_V(x)$.
\end{itemize}
We obtain $W(P,\s)\cap W_V\subset W(P_V,\s)\cap W_V$ and $W(P,\s)\cap V=\es$ for $V\in \ga$. Hence $W(P,\s)=\es$. Statement \ref{adps1} implies that $X$ is not $\dpss \cP$-nonmeager. This contradiction proves the statement.
\end{proof}

\begin{proposition} \label{pdps9}
Let $X$ be a space, $\cP\subset \sod X$ be a partible family and $\s$ be an \nfs/.
Assume that $X$ is a $\cP$-homogeneous space. Then the following conditions are equivalent:
\begin{itemize}
\item[\rm(1)]
The space $X$ is $\dpss \cP$-Baire.
\item[\rm(2)]
The space $X$ is $\dpss \cP$-nonmeager.
\end{itemize}
\end{proposition}
\begin{proof}
The implication $(1)\implies (2)$ follows from Proposition \ref{pdps2}. Let us prove $(2)\implies (1)$.
It follows from Statement \ref{adps8} that $B(\cP,\s)\neq\es$. It follows from Statement \ref{adps7} that $X$ is $\dpss \cP$-Baire.
\end{proof}

Proposition \ref{pnfs1} implies the following assertion.
\begin{proposition} \label{pdps10}
Let $X$ be a space, $\cP\subset \sod X$, $Y\subset X=\cl Y$, and $\cQ=\sorestr \cP Y$.
\begin{itemize}
\item[\rm(1)]
If $Y$ is $\dps \cQ{\s_d}$-nonmeager, then $X$ is $\dps \cP{\s_d}$-nonmeager.
\item[\rm(2)]
Let $k\in\sset{h,s}$. If $Y$ is $\dps \cQ{\s_k}$-nonmeager, then $X$ is $\dps {\pp^v(\cP)}{\s_k}$-nonmeager.
\end{itemize}
\end{proposition}

\subsection{More classes}\label{sec-dp}

\begin{definition}
Let $\s$ be a normal square functor, $X$ be a space, $\ga$ be an ordinal, and
$\cP = \sod X$. Let $\pp^k_l: \Exp{\sod X}\to \Exp{\sod X}$ be one of the mappings considered in Section \ref{sec-sonbd}:
\[
\pp^k_l\in\sset{\pp,\pp^e,\pp^v,\pp^c,\pp_+,\pp_+^v,\pp_\ga,\pp_\ga^c}.
\]
We say that $X$ is a \term{$\lrds kl{\s}$-Baire} (\term{$\lrds kl{\s}$-nonmeager}) space if $X$ is a $\dps{ \pp^k_l(\sod X)}\s$-Baire ($\dps{\pp^k_l(\sod X)}\s$-nonmeager) space. For
\[
\tD\in\sset{
\ld e,\ld v,\ld c,\rd +,\lrd c+,\rd\ga,\lrd c\ga}
\]
we have defined $\tD(\s)$-Baire ($\tD(\s)$-nonmeager) spaces. For
$k\in\sset{d,h,v,s,a}$
we say that $X$ is a \term{${\tD}_k$-Baire} (\term{${\tD}_k$-nonmeager}) space if $X$ is $\tD(\s_k)$ -Baire ($\tD(\s_k)$-nonmeager) space. For
\[
\til{\tD}\in\sset{
\ldd {e}k,\ldd {v}k,\ldd {c}k,\rdd {+}k,\lrdd c{_+} k,\rdd {_\ga}k,\lrdd c{_\ga}k}
\]
we have defined $\til{\tD}$-Baire ($\til{\tD}$-nonmeager) spaces. Also, 
by $^k\Delta^l$ without a subscript we mean $^k\Delta^l_d$:
$\lrd kl$-nonmeager ($\lrd kl$-Baire) means $\lrdd kld$-nonmeager ($\lrdd kld$-Baire). Below we give separate direct definitions for the most important classes of spaces: a space $X$ is called
\begin{itemize}
\item[\rm (1)] $\rd \ga$-nonmeager ($\rd \ga$-Baire) if $X$ is $\dps{\pp_\ga(\sod X)}{\s_d} $-nonmeager ($\dps{\pp_\ga(\sod X)}{\s_d}$-Baire);
\item[\rm (2)] $\lrd c\ga$-nonmeager ($\lrd c\ga$-Baire) if $X$ is $\dps{\pp^c_\ga(\sod X)} {\s_d}$-nonmeager ($\dps{\pp^c_\ga(\sod X)}{\s_d}$-Baire);
\end{itemize}
\end{definition}


\begin{proposition} \label{pdp1-1}
Let $\ga$ be an ordinal and $\s$ is an \nfs/. The following classes of spaces are the same:
\begin{itemize}
\item[\rm (1)] $\rd +(\s)$-nonmeager = $\rd 2(\s)$-nonmeager, $\rd +(\s)$-Baire = $\rd 2(\s)$-Baire;
\item[\rm (2)] $\lrd c\ga(\s)$-nonmeager = $\lrdd c{\ga}a$-nonmeager, $\lrd c\ga(\s)$-Baire = $\lrdd c{\ga}a$-Baire;
\item[\rm (3)] $\lrd c\ga$-nonmeager = $\lrdd c{\ga}a$-nonmeager, $\lrd c\ga$-Baire = $\lrdd c{\ga}a$-Baire.
\end{itemize}
\end{proposition}
\begin{proof} Let $X$ be a space.
Point (1) follows from the fact that $\pp_+(\sod X)=\pp_2(\sod X)$.
Point (2) follows from the fact that the elements of the family $\pp^c_\ga(\sod X)$ are closed.
Point (3) follows from (2).
\end{proof}


Propositions \ref{psonbd9},  \ref{psonbd11}, \ref{pdps9} and Statement \ref{adps8} imply the following assertion.
\begin{proposition} \label{pdp1}
Let $X$ be a space, $\s$ be an \nfs/, and $n\in\om$.
Consider the assertions:
\begin{itemize}
\item[\rm(*)]
\begin{enumerate}
\item
A space $X$ is $\tD(\s)$-nonmeager if and only if some nonempty open $U\subset X$ is $\tD(\s)$-Baire.
\item
If the space $X$ is homogeneous, then $X$ is $\tD(\s)$-nonmeager if and only if $X$ is $\tD(\s)$-Baire.
\end{enumerate}
\end{itemize}
\begin{itemize}
\item[\rm (1)] If $\tD=\rd n$, then {\rm (*)} holds.
\item[\rm (2)] If the space $X$ is quasi-regular and $\tD\in\sset{\ld v, \ld c, \lrd c+, \lrd cn }$, then {\rm (*)} holds.
\end{itemize}
\end{proposition}


\begin{statement} \label{adp2}
Let $(X,\tp)$ be a meager quasi-regular space, and let $\tps=\tp\setminus\sset\es$.
Then there are $P_n\in \sod X$ for $n\in \om$ such that
\begin{enumerate}
\item[\rm (1)] $P_n$ is nowhere dense in $X\times X$;
\item[\rm (2)] $P_0\circ P_n= P_{n}$;
\item[\rm (3)] $P_n\circ P_n\subset P_{2n}$;
\item[\rm (4)] $\cl{P_n}\subset P_{n+1}$;
\item[\rm (5)] $\bigcup\limits_{n\in\om}P_n=X\times X$.
\end{enumerate}
\end{statement}
\begin{proof}
Take $\ga_n$, $n\in\om$, as in Statement \ref{adp1}. We set
\begin{align*}
f(x)&=\max \set{m\in\om: x\in \bigcup \ga_m},
\\
f_n(x)&=\max \sset{0, f(x) - n}
\end{align*}
for $x\in X$ and $n\in \om$.
We set
\begin{align*}
P_n=\{ (x,y)\in X\times X :\
& \text{there exists a }U\in \ga_l,\text{ where }l=f_n(x),
\\
& \text{such that }x,y\in U
\}.
\end{align*}
Let us check (1). Let $U,V\in\tps$. Let $y\in V$. Let $m=f(y) + n + 1$. Take $W\in \ga_m$ such that $U'=W\cap U\neq\es$. Take $U''\in \ga_{m-n}$ such that $W\subset U''$. Then $P_n(U')\subset U''$ and $y\notin \cl{U''}$. Let $V'=V\setminus \cl{U''}$. Then $\es\neq U'\times V'\subset U\times V$ and $P_n\cap U'\times V'=\es$.

Let us check (2). Let $x\in X$, and let $m=f(x)$. If $m\leq n$, then $P_{n}(x)=X$. Consider the case $m> n$. Take $V\in\ga_{m}$ and $U\in\ga_{m-n}$ such that $x\in V \subset U$. Then $P_0(x)=V$, $P_n(x)=U$ and $P_n(V)=U$.

Let us check (3). Let $x\in X$, and let $m=f(x)$. If $m\leq 2n$, then $P_{2n}(x)=X\supset P_{n}(x)$. Consider the case $m> 2n$. Take $V\in\ga_{m-n}$ and $U\in\ga_{m-2n}$ such that $x\in V \subset U$. Then $P_n(x)=V$, $P_{2n}(x)=U$ and $P_n(V)\subset U$.

Let us check (4). By construction, $\clv{P_n}\subset P_{n+1}$. Proposition \ref{psonbd1} implies that $\cl{P_n}\subset \clv{P_0\circ P_n}$. Since, by virtue of (2), $P_n=P_0\circ P_n$, it follows that $\cl{P_n}\subset P_{n+1}$.

Let us check (5). Let $x,y\in X$, $m=f(x)$. Then $P_{2m}(x)=X$ and $(x,y)\in P_{2m}$.
\end{proof}


\begin{proposition} \label{pdp2}
Let $X$ be a quasi-regular space, $\s$ be an \nfs/ and $n\in\om$.
If $\tD\in\sset{
\rd n, \ld v, \ld c, \lrd c+, \lrd cn }$, then
\begin{enumerate}
\item[\rm(1)]
if $X$ is $\tD(\s)$-nonmeager, then $X$ is nonmeager;
\item[\rm(2)]
if $X$ is $\tD(\s)$-Baire, then $X$ is Baire.
\end{enumerate}
\end{proposition}
\begin{proof}
By Proposition \ref{pdp1} it suffices to prove (1). Assume that, on the contrary, $X$ is a meager space. Take $P_n\in\sod X$, $n\in\om$, as in Statement \ref{adp2}.
If $\cP=\sod X$ and
\[
\cQ\in\sset{\cP,\pp^v(\cP),\pp^c(\cP),\pp_+^v(\cP),\pp_n(\cP),\pp_n^c( \cP)},
\]
then $P_m\in\pp^e(\cQ)$ for some $m\in\om$. Since $P_m$ is nowhere dense in $X\times X$, $X$ is not $\dps\cQ\s$-nonmeager. This contradiction proves the proposition.
\end{proof}

\begin{proposition} \label{pdp3}
Let $X$ be a space, and let $Y\subset X=\cl Y$.
\begin{itemize}
\item[\rm(1)]
Let $\tD\in\sset{
\rd n, \ld v, \ld c, \lrd c+, \lrd cn }$.
If $Y$ is $\tD$-nonmeager, then $X$ is $\tD$-nonmeager.
\item[\rm(2)]
Suppose that each point of the set $Y$ is a point of semi-regularity of the space $X$ (for example, $X$ is a semi-regular space).
If $Y$ is $\dd s$-nonmeager, then $X$ is $\dd s$-nonmeager.
\end{itemize}
\end{proposition}
\begin{proof}
Assertion (1) follows from Propositions \ref{psonbd4+1} and \ref{pdps10}. Let us prove (2).
Let $P\in\sod X$. Take $Q\in\sod X$ such that ${Q(y)}\subset P(y)$ and $Q(y)$ is a regular open set for $y\in Y$. Since $Y$ is a $\dd s$-nonmeager space, there exists a nonempty open $U\subset X$ such that $R\in \s_s(S)$, where $S=Y\cap U$ and $R= \sorestr QS$.
\begin{lemma*} Let $y\in S$. If $S=R(y)$, then $U\subset Q(y)\subset P(y)$.
\end{lemma*}
\begin{proof}
$U\subset \Int{\cl S}\subset \Int{\cl {R(y)}}\subset \Int{\cl {Q(y)}}=Q(y)$.
\end{proof}
The lemma, the relation $R\in \s_s(S)$, and Proposition \ref{pnfs1} imply $\sorestr PU\in \s_s(U)$.
\end{proof}

Propositions \ref{pdp1}, \ref{pdp2}, and \ref{pdp3}  imply the following assertion.

\begin{theorem} \label{tdp1}
Let $X$ be a quasi-regular space and $n\in\om$.
If 
$\tD\in\sset{\rd n, \lrd cn}$, 
then
\begin{enumerate}
\item[\rm(1)]
if $X$ is $\tD$-nonmeager ($\tD$-Baire), then $X$ is nonmeager (Baire);
\item[\rm(2)]
a space $X$ is $\tD$-nonmeager if and only if some nonempty open $U\subset X$ is $\tD$-Baire;
\item[\rm(3)]
a space $X$ is $\tD$-Baire if and only if any nonempty open $U\subset X$ is $\tD$-nonmeager;
\item[\rm(4)]
if the space $X$ is homogeneous, then $X$ is $\tD$-nonmeager if and only if $X$ is $\tD$-Baire.
\item[\rm(5)]
if $Y\subset X=\cl Y$ and $Y$ is $\tD$-Baire,  then $X$ is $\tD$-Baire.
\end{enumerate}
\end{theorem}

\section{Modifications of the Banach-Mazur game}\label{sec-bmg}

In this section we use topological games. Basic concepts can be found in \cite{Oxtoby1957,tel-gal-1986,arh2010,rezn2008,rezn2022-2008}.
The terminology and notation correspond to \cite{rezn2022gbtg}.

Let $G$ be a game in which player $\al$ has a winning strategy.
Let us call this game
 $G$ \term{$\al$-favorable}. If there is no such strategy, then  $G$ is called \term{$\al$-unfavorable}.

If the definition of a game $G$ depends on only one parameter, namely, some space $X$, that is, the game $G=\ggg(X)$, then we write \term{space $X$ $(\al,\ggg)$-favorable} if the game $\ggg(X)$ is $\al$-favorable and \term{space $X$ is $(\al,\ggg)$-unfavorable} if the game $\ggg(X)$   is $\al$-unfavorable.

We need games from \cite{rezn2022gbtg}, we will repeat the definition of those games that are directly used in the proofs. The rest of the games can be found in \cite{rezn2022gbtg}.

Let $(X,\tp)$ be a space, $\tps=\tp\setminus\sset\es$, $U\in\tps$.
We denote
\begin{align*}
\vv(X)\eqdef \{\sqnn V\in \tps^\om:\ & V_{n+1}\subset V_n\text{ for }\nom\},
\\
\uu(U)\eqdef\{V\in\tps:\ & V\subset U\}.
\end{align*}

\paragraph{The games $BM(X,\cV)$ and $MB(X,\cV)$ \cite{rezn2022gbtg}}
Let $\cV\subset \vv(X)$.
There are two players, $\al$ and $\be$. 
These games differ in the first move of player $\al$.
On the first move, player $\al$ chooses $U_0=X$ in the game $BM(X,\cV)$ and $U_0\in\uu(X)$ in the game $MB(X,\cV)$. Player $\be$ chooses $V_0\in\uu(U_0)$. On the $n$th move, $\al$ chooses $U_n\in\uu(V_{n-1})$ and $\be$ chooses $V_n\in\uu(U_n)$. After a countable number of moves, the winner is determined: player $\al$ wins if $\sqnn V\in\cV$.


\paragraph{The games $OD_l(X,\pn)$ and $DO_l(X,\pn)$}
Let $\pn$ be a $\pi$-net of $X$.
There are two players, $\al$ and $\be$. 
These games differ in the first move of player $\al$.
On the first move, player $\al$ chooses $U_0=X$ in the game $OD_l(X,\pn)$ and $U_0\in\uu(X)$ in the game $DO_l(X,\pn)$. Player $\be$ chooses $V_0\in\uu(U_0)$ and $M_0\in\pn$, $M_0\subset U_0$. On the $n$th move $\al$ chooses $U_n\in\uu(V_{n-1})$ and $\be$ chooses $V_n\in\uu(U_n)$ and $M_n\in\pn $, $M_n\subset U_n$. After a countable number of moves, the winner is determined: player $\al$ wins if $\ltu_{\nom} M_n \cap \bigcap_{\nom} V_n\neq\es$.


We denote
\begin{align*}
\vvbm(X)\eqdef \{\sqnn V\in \vv(X):\ & \bigcap_{\nom} V_n\neq\es \},
\\
\vv_f(X)\eqdef \{\sqnn V\in\vv (X) :\
& \text{for some }x\in \bigcap_{\nom} V_n \text { the family }
\\
&\sqnn V\text { forms a base at the point }x \},
\\
\pn_o(X)\eqdef \tps,\quad & \pn_p(X)\eqdef \set{\sset x:x\in X}.
\end{align*}
We define the games
\begin{align*}
BM(X) &\eqdef BM(X,\vvbm(X)) & MB(X) &\eqdef MB(X,\vvbm(X)),
\\
BM_f(X) &\eqdef BM(X,\vv_f(X)) & MB_f(X) &\eqdef MB(X,\vv_f(X)),
\\
OD_{o,l}(X) &\eqdef OD_l(X,\pn_o(X)) & DO_{o,l}(X) &\eqdef DO_l(X,\pn_o(X)),
\\
OD_{p,l}(X) &\eqdef OD_l(X,\pn_p(X)) & DO_{p,l}(X) &\eqdef DO_l(X,\pn_p(X)).
\end{align*}

The game $BM(X)$ is the classical Banach--Mazur game.

\begin{theorem}[Banach--Oxtoby \cite{Oxtoby1957}; see also \cite{tel-gal-1986,arh2010}] \label{tbmb-bm}
Let $X$ be a space.
\begin{itemize}
\item[\rm(1)]
$X$ is Baire if and only if $BM(X)$ is $\be$-unfavorable;
\item[\rm(2)]
$X$ is nonmeager if and only if $MB(X)$ is $\be$-unfavorable.
\end{itemize}
\end{theorem}

\begin{definition}
We define classes of spaces:
\begin{align*}
\text{$\gbm_f$-Baire} & \eqdef \text{$(\be,BM_f)$-unfavorable;}
&
\text{$\gbm_f$-nonmeager} & \eqdef \text{$(\be,MB_f)$-unfavorable;}
\\
\text{$\god_{o,l}$-Baire} & \eqdef \text{$(\be,OD_{o,l})$-unfavorable;}
&
\text{$\god_{o,l}$-nonmeager} & \eqdef \text{$(\be,DO_{o,l})$-unfavorable;}
\\
\text{$\god_{p,l}$-Baire} & \eqdef \text{$(\be,OD_{p,l})$-unfavorable;}
&
\text{$\god_{p,l}$-nonmeager} & \eqdef \text{$(\be,DO_{p,l})$-unfavorable.}
\end{align*}
\end{definition}

\subsection{Strategies vs tactics}\label{sec-gtt}
Let $\ggg$ be a topological game played by players $\al$ and $\be$. In a typical topological game, the game $\ggg(X)$ is a \term{sequential infinite game of perfect information}. The game consists of players taking turns making moves. After a countable number of moves, after analyzing which players made moves, the winner is determined. When making a move, a player makes a choice of move based on all previous moves of both players.
These choices determine the player's \term{strategy}.

A strategy in a sequential game is called a \term{tactic} \cite{Choquet1969} if the player takes into account only the opponent's previous move in his choice. Such strategies are also called $1$-tactics or stationary strategies.
A strategy is called a \term{limited information strategy} if the player takes into account not all, but only some previous moves of the opponent and himself.
Thus, a tactic is a limited information strategy. In addition to tactics, other limited information strategies are also considered \cite{tel-gal-1986,BrianAlanMilovich2021}.

Let $\cP$ and $\cQ$ is some topological properties. An important role when using topological games is played by theorems of the form:
\begin{itemize}
\item[($\ggg_1$)] if $X$ holds $\cP$, then $\be$ has a winning strategy in the game $\ggg(X)$.
\end{itemize}
A typical proof of theorems of the form ($G_1$) is that the $\cP$ property is used to construct a strategy for $\be$. Often, when building a strategy, not all previous moves are used, but only a certain part.
That is, for $\be$ there is some limited information strategy. Sometimes a winning tactic is built for $\be$. In the next section, we will find winning tactics for player $\be$ in the games we are interested in.

Theorem $(\ggg_1)$ is equivalent to
\begin{itemize}
\item[($\ggg_2$)] if $X$ holds $\cP$, $X$ is $(\ggg,\be)$-favorable.
\end{itemize}
If $\cQ=\lnot \cP$, then the theorem ($\ggg_1$) is equivalent to
\begin{itemize}
\item[($\ggg_3$)] if $X$ is $(\ggg,\be)$-unfavorable, then $\cQ$ holds for $X$.
\end{itemize}
We use the form ($\ggg_3$), since for the games considered in this article, the $(\ggg,\be)$-unfavorability property describes nonmeager-type and Baire-type properties.

Let us denote by $\tggg(X)$ a game in which the rules are the same as in $\ggg(X)$, but the player $\be$ uses only tactics in the game. The game $\tggg(X)$ is a sequential infinite game of imperfect information.

\begin{theorem}\label{t:svst:1}
Let $\ggg(X)$ be a topological game. If $X$ is $(\ggg,\be)$-unfavorable, then $X$ is $(\tggg,\be)$-unfavorable.
\end{theorem}
\begin{proof}
Let $s_\be$ be some tactics of player $\be$. Since a tactic is a strategy and $X$ is $(\ggg,\be)$-unfavorable, there is a strategy $s_\al$ of player $\al$ such that after playing using strategies $s_\al$ and $s_ \be$ player $\be$ will lose.
\end{proof}

In Section \ref{sec-gd} we will prove theorems of the form
\begin{itemize}
\item[($\tggg_3$)] if $X$ is $(\tggg,\be)$-unfavorable, then $\cQ$ holds for $X$.
\end{itemize}

From Theorem \ref{t:svst:1} it follows that the theorem ($\tggg_3$) is stronger than the theorem ($\ggg_3$).
The next question arises.
\begin{itemize}
\item[($Q_1$)] Are the topological properties of $X$ is $(\ggg,\be)$-unfavorable and $X$ is $(\tggg,\be)$-unfavorable the same?
\end{itemize}
This question is equivalent to the following question.
\begin{itemize}
\item[($Q_2$)] Suppose that player $\be$ has a winning strategy in the game $\ggg(X)$. Does $\be$ have a winning tactic in this game?
\end{itemize}

The question ($Q_2$) has been studied for many years. The article \cite{tel-gal-1986} is devoted to this problem. In \cite{tel-gal-1986} the general Theorem 1 was proved, which allows for a number of games to solve ($Q_2$) positively, for example, for the Banach-Mazur game. But for many games this issue is resolved negatively \cite{tel-gal-1986,BrianAlanMilovich2021,Debs1985}.

\subsection{Playing against tactics}\label{sec-gtt}


Let $(X,\tp)$ be a space and $\tps=\tp\setminus\sset\es$, $\pn$ be a $\pi$-net $X$.
We denote
\begin{align*}
T(X,\pn)&\eqdef \set{\psi\in\pn^\tps: M\subset U \text{ for } U\in\tps\text{ and }M=\psi(U)} ,
\\
\tbm X&\eqdef T(X,\tps),
\\
\tod X\pn&\eqdef T(X,\tps) \times T(X,\pn).
\end{align*}

Each $\ph\in \tbm X$ uniquely corresponds to the tactics of player $\be$ in the games $BM(X,\cV)$ and $MB(X,\cV)$: on the $n$-th move player $\be$ selects $V_n=\ph(U_n)$. Each $(\ph,\psi)\in \tod X\pn$ uniquely corresponds to the tactics of player $\be$ in the games $OD_l(X,\pn)$ and $DO_l(X,\pn)$: on the $n$th move, player $\be$ chooses $V_n=\ph(U_n)$ and $M_n=\psi(U_n)$.

The modification of the games $BM$, $MB$, $OD_l$ and $DO_l$, in which the player $\be$ uses only tactics, will be denoted by $\tg BM$, $\tg MB$, $\tg OD_l$, $ \tg DO_l$, respectively.

We will assume that
\begin{itemize}
\item
$\tbm X$ is the set of strategies of player $\be$ in the games $\tg BM$ and $\tg MB$;
\item
$\tod X\pn$ is the set of strategies of player $\be$ in the games $\tg OD_l$ and $\tg DO_l$.
\end{itemize}

We will give variants of definitions in which games are replaced by games ``with a hat'', that is, those in which the player $\be$ uses only tactics.
We have games
\begin{align*}
\tg BM_f, && \tg MB_f, && \tg OD_{o,l}, && \tg DO_{o,l},
&& \tg OD_{p,l}, && \tg DO_{p,l}
\end{align*}
and classes of spaces:
\begin{align*}
\text{$\gbmt_f$-Baire;} &&
\text{$\gbmt_f$-nonmeager;} &&
\text{$\godt_{o,l}$-Baire;} &&
\text{$\godt_{o,l}$-nonmeager;} &
\\
&&
\text{$\godt_{p,l}$-Baire;} &&
\text{$\godt_{p,l}$-nonmeager.}
\end{align*}

Theorem \ref{t:svst:1}  implies the following assertion.
\begin{proposition}\label{pgtt4}
Let $X$ be a space.
Then
\begin{itemize}
\item if $X$ is {$\gbm_f$-nonmeager} ($\gbm_f$-Baire), then $X$ is {$\gbmt_f$-nonmeager} ($\gbmt_f$-Baire);
\item if $X$ is {$\god_{o,l}$-nonmeager} ($\god_{o,l}$-Baire), then $X$ is {$\godt_{o,l}$-nonmeager} ($\godt_{o,l}$-Baire);
\item if $X$ is {$\god_{p,l}$-nonmeager} ($\god_{p,l}$-Baire), then $X$ is {$\godt_{p,l}$-nonmeager} ($\godt_{p,l}$-Baire).
\end{itemize}
\end{proposition}

The following proposition is easy to deduce from the definitions.

\begin{proposition}\label{pgtt3}
Let $X$ be a space.
\begin{itemize}
\item
Let $\cV\subset \vv(X)$.
A game $\tg MB(X,\cV)$ ($\tg BM(X,\cV)$) is $\be$-unfavorable
 if and only if
for any $\ph\in\tbm X$ there exists $\sqn{U_n}\in \vv(X)$ (for which $U_0=X$) such that
\begin{itemize}
\item
$V_n=\ph(U_n)$, $U_{n+1}\subset V_n$ for $\nom$;
\item
$\sqn{V_n}\in \cV$.
\end{itemize}
\item
Let $\pn$ be a $\pi$-net $X$.
A game $\tg DO(X,\pn)$ ($\tg OD(X,\pn)$) is $\be$-unfavorable  if and only if
for any $(\ph,\psi)\in\tod X\pn$ there exists $\sqn{U_n}\in \vv(X)$ (for which $U_0=X$) such that
\begin{itemize}
\item
$V_n=\ph(U_n)$, $M_n=\psi(U_n)$ and $U_{n+1}\subset V_n$ for $\nom$;
\item
$\ltu_{\nom} M_n \cap \bigcap_{\nom} V_n\neq\es$.
\end{itemize}
\end{itemize}
\end{proposition}

\section{Relationship between $\D$ and $\ggg$ Baire}\label{sec-gd}

In this section, a connection is established between Baire-type and nonmeager-type properties, determined by using the diagonal and topological games.

\begin{proposition}\label{pgd1} Let $(X,\tp)$ be a space, and let $\tps=\tp\setminus\sset\es$.
\begin{itemize}
\item[{\rm (1)}]
If $X$ is $\godt_{o,l}$-nonmeager ($\godt_{o,l}$-Baire), then $X$ is $\D$-nonmeager ($\D$-Baire).
\item[{\rm (2)}]
If $X$ is $\godt_{p,l}$-nonmeager ($\godt_{p,l}$-Baire), then $X$ is $\D_h$-nonmeager ($\D_h$-Baire).
\item[{\rm (3)}]
If $X$ is $\gbmt_{f}$-nonmeager ($\gbmt_{f}$-Baire), then $X$ is $\D_s$-nonmeager ($\D_s$-Baire).
\end{itemize}
\end{proposition}
\begin{proof} We will prove the proposition by contradiction based on Proposition \ref{pgtt3}. We carry out the proof for strengthening of nonmeager property; for strengthening of the Baire property, the proof is the same  with the only difference that we must additionally ensure during the construction that $V_0=\ph(X)$ is not strengthened nonmeager.

(1)
Assume that $X$ is not $\D$-nonmeager. Then there is an semi-neighborhood of the diagonal $P\in \sod X$ such that $\sorestr PU\notin \s_d(U)$ for any $U\in\tps$. 
That is, $U\times U\not\subset \cl P$ for any $U\in\tps$.
We define $(\ph,\psi)\in\tod X{\pn_o(X)}$. Recall that $\pn_o(X)=\tps$. Let $U\in\tps$. 
Since 
$U\times U\not\subset \cl P$
it follows  by Proposition \ref{pnfs1} that there exists $V\subset U$, $V\in\tps$, such that $U\not\subset \cl{P(V)}$. Let $\ph(U)=V$ and $\psi(U)=M$, where $M=U\setminus \cl{P(V)}$. Note that $M \cap P(V)=\es$. We have constructed mappings $\ph$ and $\psi$.
By Proposition \ref{pgtt3}, there exists $\sqn{U_n}\in \vv(X)$ such that 
\begin{itemize}
\item[(a)]
$V_n=\ph(U_n)$, $M_n=\psi(U_n)$, and $U_{n+1}\subset V_n$ for $\nom$;
\item[(b)]
$\sqn{V_n,M_n}\in \ww_l(X)$, i.e., $\ltu_{\nom} M_n \cap \bigcap_{\nom} V_n\neq\es$.
\end{itemize}
The construction of $\psi$ implies
\begin{itemize}
\item[(c)]
$M_n \cap P(V_n)=\es$ for any $\nom$.
\end{itemize}
Let $x_*\in \ltu_{\nom} M_n \cap \bigcap_{\nom} V_n$. Then $M_n \cap P(x_*) \ne \es$ for some $\nom$. We obtain $M_n \cap P(V_n)\neq\es$. This contradictions (c).

(2) 
Assume that $X$ is not $\D_h$-nonmeager. Then there is an almost neighborhood of the diagonal $P\in \sod X$ such that $\sorestr PU\notin \s_h(U)$ for any $U\in\tps$. 
That is, $U\times U\not\subset \clh P$ for any $U\in\tps$.
We define $(\ph,\psi)\in\tod X{\pn_p(X)}$. Recall that $\pn_p(X)=\set{\sset x:x\in X}$. Let $U\in\tps$. Since 
$(U\times U)\cap P= \sorestr PU\notin \s_h(U)$,
it follows  by Proposition \ref{pnfs1}, there exists $V\subset U$, $V\in\tps$, such that $U\not\subset { P(V)}$. Let $\ph(U)=V$ and $\psi(U)=\sset x$, where $x\in U\setminus {P(V)}$. Note that $x\notin P(V)$. We have constructed mappings $\ph$ and $\psi$.
By Proposition \ref{pgtt3}, there exists $\sqn{U_n}\in \vv(X)$ such that 
\begin{itemize}
\item[(a)]
$V_n=\ph(U_n)$, $\sset{x_n}=\psi(U_n)$, and $U_{n+1}\subset V_n$ for $\nom$;
\item[(b)]
$\sqn{V_n,\sset{x_n}}\in \ww_l(X)$, i.e., $\ltu_{\nom} x_n \cap \bigcap_{\nom} V_n\neq\es$.
\end{itemize}
The construction of $\psi$ implies
\begin{itemize}
\item[(c)]
$x_n \notin P(V_n)$ for any $\nom$.
\end{itemize}
Let $x_*\in \ltu_{\nom} x_n \cap \bigcap_{\nom} V_n$. Then $x_n \in P(x_*)$ for some $\nom$. We obtain $x_n \in P(V_n)$. This contradicts (c).

(3)
Assume that $X$ is not $\D_s$-nonmeager. Then there is an almost neighborhood of the diagonal $P\in \sod X$ such that $(U\times U)\cap P=\sorestr PU\notin \s_h(U)$ for any $U\in\tps$. Let us define $\ph\in\tbm X$. Let $U\in\tps$. Since $\sorestr PU\notin \s_s(U)$, by Proposition \ref{pnfs1}, there exists $V\subset U$, $V\in\tps$, such that
$U\not\subset P(x)$ for any $x\in V$.
 Let $\ph(U)=V$. We have constructed a mapping $\ph: \tps\to\tps$.
By Proposition \ref{pgtt3}, there exists $\sqn{U_n}\in \vv(X)$ such that 
\begin{itemize}
\item[(a)]
$V_n=\ph(U_n)$, $U_{n+1}\subset V_n$ for $\nom$;
\item[(b)]
$\sqnn V \in \vv_f(X)$, i.e., $\bigcap_{\nom} V_n\neq\es$ and the family $\sqnn V$ is a base of some point $x_*\in \bigcap_{\nom} V_n$.
\end{itemize}
The construction of $\ph$ implies
\begin{itemize}
\item[(c)]
$U_n\not\subset P(x)$ for any $\nom$ and $x\in V_n$.
\end{itemize}
Note that the family $\sqnn U$ is also a base of the point $x_*$. There exists $\nom$ such that $U_n\subset P(x_*)$. Since $x_*\in V_n$, we obtain a contradiction to (c).
\end{proof}

Propositions \ref{pgd1}, \ref{pdps7}, \ref{pgtt4}, 
\cite[Proposition 15]{rezn2022gbtg}
and \cite[Theorems 3 and 4]{rezn2022gbtg}
imply the following theorem.
\begin{theorem} \label{tgd1}
Let $X$ be a space.
In the diagram below, the arrow from $A$ to $B$
means that
\begin{itemize}
\item[\rm (1)] If $X$ is an $A$-nonmeager space, then $X$ is a $B$-nonmeager space;
\item[\rm (2)] if $X$ is an $A$-Baire space, then $X$ is a $B$-Baire space.
\end{itemize}

{
\def\a#1{\arrow[#1]}
\[
\begin{tikzcd}
\gbm_{f} \a r \a d & \gbmt_{f} \a r \a d& \D_s \a d\\
\god_{p,l} \a r \a d& \godt_{p,l} \a r \a d& \D_h \a d\\
\god_{o,l} \a r & \godt_{o,l} \a r & \D \\
\end{tikzcd}
\]}

Let $\cP_k$ ($\cP_c$) be the smallest class of spaces that
\begin{itemize}
\item contains $p$-spaces and strongly $\Sigma$-spaces;
\item is closed under arbitrary (countable) products;
\item  is closed under taking open subspaces.
\end{itemize}

Given $\wt\D\in \sset{\D_s,\D_h,\D}$, if a regular Baire (nonmeager) space $X$ belongs to the class of spaces described in $(\wt\D)$, then $ X$ is $\wt\D$-Baire ($\wt\D$-nonmeager).
\begin{itemize}
\item[$(\D_s)$] $\sigma$-spaces.
\item[$(\D_h)$] is the smallest class of spaces that
\begin{itemize}
\item contains $\Sigma$-spaces and $w\D$-spaces;
\item is closed under products by spaces from the class $\cP_c$;
\item  is closed under taking open subspaces.
\end{itemize}
\item[$(\D)$] is the smallest class of spaces that
\item contains $\Sigma$-spaces, $w\D$-spaces, and feebly compact spaces;
\item is closed under products by spaces from the class $\cP_k$;
\item  is closed under taking open subspaces.
\end{itemize}

\end{theorem}

\section{\cdp/-Baire spaces}\label{sec-cdp}

Let $X$ be a space, and let $\cG$ be a сollection of families of open subsets of $X$.
We denote
\begin{align*}
B(X,\cG) \eqdef \{ x\in X:\ & \text{ the family }\\
&\set{\st(x,\ga):\ga\in\cG\text{ and }x\in\bigcup\ga}
\\
&\ \text{ is a base at  }x \}
\\
\dev X\eqdef \min \{\, |\cG| :\
&\cG \text{ is a сollection of open covers of }X
\\
& \text{ and }X=B(X,\cG)
,\}.
\end{align*}

Let $Y\subset X$, and let $\cN$ be the family of all nowhere dense subsets of $X$. If $Y\not \subset \bigcup \cN$, then we put $\nov {Y,X} \eqdef \infty$, otherwise we put
\begin{align*}
\nov {Y,X} &\eqdef \min \set{|\cL|: \cL\subset \cN\text{ and } Y\subset \bigcup\cL},
\\
\nov {X} &\eqdef \nov {X,X}.
\end{align*}

\begin{definition}
Let $X$ be a space.
\begin{itemize}
\item
We call a space a \term{\cdp/-space} if $\dev X<\nov X$, that is, if there exists a сollection $\cP$ of open coverings of the space $X$ such that $|\cP|<\nov X $ and $B(X,\cP)=X$.
\item
We call a space a \term{\cdpz/-space} if there exists a сollection $\cP$ of open partitions of $X$ such that $|\cP|<\nov X$ and $B(X,\cP)=X$ .
\item We call a space $X$ \term{\cdp/-nonmeager} if there exists a сollection $\cP$ of open families in $X$ such that $|\cP|<\nov{Y,X}$ for $Y =B(X,\cP)$.
\item We call a space $X$ \term{\cdp/-Baire} if
every nonempty open subset of $X$ is a \cdp/-nonmeager space.
\end{itemize}
\end{definition}

The authors of \cite{cdp2010} studied spaces for which $\dev X<\nov X$, that is, \cdp/-spaces. Metrizable nonmeager spaces are \cdp/-spaces. Examples of non-metrizable \cdp/-spaces can be obtained using Martin's axiom ($MA$).

\begin{proposition}[Theorem 2.3 \cite{tall1974}] \label{pcdp1}
($MA$)
Let $\tau<2^\om$ be an infinite cardinal, and let $X$ be an absolute $G_\tau$ space with countable Souslin number. Then $\nov X\geq 2^\om$.
\end{proposition}

Proposition \ref{pcdp1}  implies the following assertion.
\begin{proposition} \label{pcdp2}
($MA$)
Let $\tau<2^\om$ be an infinite cardinal, and let $X$ be an absolute $G_\tau$ space with countable Souslin number and $\dev X\leq \tau$, for example $X=\R^\tau$. Then $X$ is a \cdp/-space.
\end{proposition}

Clearly, a \cdp/-space is a \cdpz/-space.

\begin{proposition} \label{pcdp3}
Let $X$ be a \cdp/-nonmeager space.
\begin{itemize}
\item[\rm(1)]
Then some nonempty open subset $U\subset X$ contains a dense \cdpz/-space $Y\subset U\subset \cl Y$.
\item[\rm (2)] If $X$ is a $\pi$-semiregular space, then $U$ and $Y$ can be chosen so that $X$ is semiregular at every point in $Y$.
\end{itemize}
\end{proposition}
\begin{proof}
Let $\cG=\set{\ga_\al:\al<\tau}$ be a family of open families in $X$ such that $\tau=|\cG|<\nov{S,X}$ for $S=B(X,\cG)$.

Let $\al<\tau$. There is an open disjoint family $\nu_\al$ such that
\begin{itemize}
\item
if $U\in \nu_\al$ and $U\cap\bigcup \ga_\al\neq \es$, then $U\subset V$ for some $V\in \ga_\al$;
\item
$\cl{\bigcup\nu_\al}=X$.
\end{itemize}
If $X$ is a $\pi$-semiregular space, then one can additionally ensure the condition
\begin{itemize}
\item
the family $\nu_\al$ consists of regular open sets.
\end{itemize}
Let $G_\al=\bigcup\nu_\al$ and $F_\al=X\setminus G_\al$. The set $F_\al$ is nowhere dense in $X$.

Let $\cV=\set{\nu_\al: \al<\tau}$, $G=\bigcap\set{G_\al:\al<\tau}$, $F=\bigcup\set{F_ \al:\al<\tau}$,
\[
Z=G\cap B(X,\cV).
\]
Then $S\setminus F\subset Z$ and $\nov{Z,X}\geq \nov{Y,X}>\tau$. Let $U=\Int{\cl Z}$. Since $\nov{Z,X}$ is uncountable, $U\neq\es$. Let $Y=Z\cap U$.
The space $Y$ is a \cdpz/-space, $Y\subset U\subset \cl Y$. If the families $\nu_\al$ consisted of regular open sets, then $X$ is semiregular at every point in $Y$.
\end{proof}

\begin{proposition} \label{pcdp4}
Let $X$ be a space, $U$ be an open subset of $X$, and $Y\subset U\subset \cl Y$. If $Y$ is \cdp/-Baire and $X$ is semiregular at points in $Y$, then $X$ is \cdp/-nonmeager.
\end{proposition}
\begin{proof}
Take a сollection $\cP$ of open in $Y$ coverings of the space $Y$ such that $|\cP|<\nov Y$ and $B(Y,\cP)=Y$. For $\ga\in\cP$ we set $\til \ga = \set{\Int{\cl V}: V\in\ga}$.
We put $\wt \cP=\set{\til \ga: \ga\in\cP}$.
Then $Y\subset \wt Y = B(X,\wt\cP)$ and $|\wt\cP| < \nov{Y} \leq \nov{\wt Y,X}$.
\end{proof}

Propositions \ref{pcdp3} and \ref{pcdp4} imply  the following theorem.

\begin{theorem} \label{tcdp1}
Let $X$ be a semiregular space. The following conditions are equivalent:
\begin{itemize}
\item[\rm (1)] $X$ is a \cdp/-nonmeager space;
\item[\rm (2)] some nonempty open subset $U\subset X$ contains a dense \cdp/-space $Y\subset U\subset \cl Y$;
\item[\rm (3)] some nonempty open subset $U\subset X$ contains a dense \cdpz/-space $Y\subset U\subset \cl Y$.
\end{itemize}
The following conditions are equivalent:
\begin{itemize}
\item[\rm (4)] $X$ is a \cdp/-Baire space;
\item[\rm (5)]  $X$ contains a dense \cdp/-space $Y\subset X\subset \cl Y$;
\item[\rm (6)] $X$ contains a dense \cdpz/-space $Y\subset X\subset \cl Y$.
\end{itemize}
\end{theorem}

Theorem \ref{tcdp1} implies

\begin{proposition} \label{pcdp5}
Let $X$ be a semiregular space.
If $X$ has a metrizable nonmeager subspace, then $X$ is a \cdp/-nonmeager space.
\end{proposition}


Theorem \ref{tgd1} implies that a $\gbm_{f}$-Baire ($\gbm_{f}$-nonmeager) space is a $\D_s$-Baire ($\D_s$-nonmeager) space. 
Proposition 15 and Theorem 3 from \cite{rezn2022gbtg}
give the following examples of $\gbm_{f}$-Baire ($\gbm_{f}$-nonmeager) spaces; Baire (nonmeager) spaces from the classes of spaces listed below are $\gbm_{f}$-Baire ($\gbm_{f}$-nonmeager):
\begin{itemize}
\item[]
metrizable spaces, Moore spaces, developable space, semiregular $\sigma$-spaces, and semiregular spaces with a countable network.
\end{itemize}

\begin{proposition}\label{pdsp1}
Let $X$ be a space. If $X$ is \cdp/-nonmeager (\cdp/-Baire), then $X$ is $\D_s$-nonmeager ($\D_s$-Baire).
\end{proposition}
\begin{proof}
Let us prove the propositions for the case when $X$ is \cdp/-nonmeager. Let $\cG$ be a family of open families in $X$ such that $|\cG|<\nov{Y,X}$ for $Y=B(X,\cG)$.
Let $P\in \sod X$ be an almost neighborhood of the diagonal $X\times X$. For every $y\in Y$ there is $\ga_y\in\cG$ such that $\st(y,\ga_y)\subset P(y)$. We put $Y_\ga=\set{y\in Y: \ga_y=\ga}$. Since $|\cG|<\nov{Y,X}$ , it follows that $Y_{\ga}'$ is nowhere dense for some $\ga'\in \cG$. We choose $y'\in V\cap Y$, where $V=\Int(\cl{Y_{\ga'}})$. Let $W \in \St(y',\ga')$. Let $U=W\cap V$ and $M=U\cap Y$. Then $M$ is dense in $U$ and for $y\in M$
\[
U\subset W \subset \st(y,\ga')\subset P(y).
\]
Hence $M\times U\subset P$.
\end{proof}

Examples of $\D_s$-nonmeager spaces that are obtained from $\gbm_{f}$-nonmeager spaces always contain points with the first countability axiom. Propositions \ref{pdsp1} and \ref{pcdp2} imply that under the assumption MA+$\lnot$CH the space $\R^{\om_1}$ is $\D_s$-nonmeager.
\section{Continuity of operations in groups}\label{sec-tg}

The main applications of the considered generalizations of Baire spaces relate to the topic of continuity of operations in groups with topology. Such applications were considered in \cite{Reznichenko2024rtg}. In this section, we formulate results that will help to distinguish between different types of Baire property.

Let $G$ be a topologized group. A group $G$ is called \term{semitopological} if multiplication in $G$ is separately continuous, that is, for $g\in G$ the maps $x \mapsto gx$ and $x \mapsto xg$ are continuous. A group $G$ is called \term{paratopological} if the multiplication in $G$ is continuous. A group $G$ is called \term{quasitopological} if $G$ is a semitopological group and the operation of taking the inverse of $x \mapsto x^{-1}$ is continuous. A group $G$ is called \term{topological} if the operations in $G$ are continuous, that is, $G$ is a paratopological quasi-topological group.

Let $\cN_e$ be a family of open neighborhoods of a identity $e$ of $G$, $\cN(U)=\set{ (g,h)\in G\times G: h\in Ug }$ for $U\subset G$. Then $P(g)=Ug$ for $g\in G$. 
If $G$ is semitopological and $U\in \cN_e$ then $P$ is a semi-neighborhood of the diagonal.

\begin{statement} \label{s:tg:1}
Let $G$ is a semitopological group, $V,U\in\cN_e$. Then
\rm(1) $\cN(\cl U)=\clv{\cN(U)}$;
\rm(2) $\cN(UV)= \cN(V) \circ \cN(U)$;
\rm(3) $\cl{\cN(U)}\subset \cN(\cl{UV})$;
\rm(4) if $\cl{V^3}\subset U$ then $\cl{\cN(V) \circ\cN(V)}\subset \cN(U)$.
\end{statement}
\begin{proof}
(1) Obviously.

(2) $(x,y)\in \cN(UV)$ $\lrarr$ $y=uvx$ for some $u\in U$ and $v\in V$ $\lrarr$ $z=vx$ and $y=uz$ for some $u\in U$ and $v\in V$ $\lrarr$ $(x,z)\in \cN(V)$ and $(z,y)\in \cN(U)$ for some $z\in G$ $\lrarr$ $(x,y)\in \cN(V) \circ \cN(U)$.

(3) From (1) and (2) it follows $\clv{\cN(V) \circ \cN(U)} = \cN(\cl{UV})$. Proposition \ref{psonbd1} implies $\cl{\cN(U)}\subset \clv{\cN(V) \circ \cN(U)}$. Therefore, $\cl{\cN(U)}\subset \cN(\cl{UV})$.

(4) From (2) and (3) it follows $\cl{\cN(V) \circ\cN(V)}=\cl{\cN(V^2)}\subset \cN(\cl{V^3})\subset \cN(U)$.
\end{proof}

\begin{theorem}\label{ttg1}
Let $G$ be a semiregular topologized group. 
Assume that one of the following conditions is met:
\begin{itemize}
\item[\rm(1)]
$G$ is a $\D$-Baire paratopological group;
\item[\rm(2)]
$G$ is a $\D_h$-Baire quasitopological group.
\end{itemize}
Then $G$ is a topological group. 
\end{theorem}
\begin{proof}
Let $O\in \cN_e$ be a regular open set.

(1) We need to find $V\in\cN_e$ for which $V^{-1}\subset O$. Take $U\in\cN_e$ such that $U^2\subset O$. 
Then $\cl{UU}\subset \cl{O}$ and 
Statement \ref{s:tg:1} implies $\cl{\cN(U)}\subset \cN(\cl{O})$. Since $G$ is $\D$-nonmeager, there is an open nonempty $W\subset G$ such that $W\times W \subset \cN(\cl{O})$. Then $V=WW^{-1}\subset \cl{O}$. 
Since  $O$ is regular open and $V$ is open then $V\subset\Int{\cl{O}}=O$ and $V=V^{-1}\subset O$.

(2) We need to find a $V\in\cN_e$ for which $VV\subset O$. Take $U\in\cN_e$ such that $U=U^{-1}\subset O$.
 Then $P=\cN(U)$ is a semi-neighborhood of the diagonal. Since $G$ is $\D_h$-nonmeager and $P=P^{-1}$, there is an open nonempty $W\subset G$ such that $W \subset \cl{P(x)}$ for $x\in W$. Hence $WW^{-1}\subset \cl U$. Let $V=V^{-1}\in\cN_e$  such that $Vg\subset W$ for some $g\in W$. Then $VV\subset \cl U\subset O$. Since  $O$ is regular open then $VV\subset O$.
%
\end{proof}

Theorem \ref{ttg1} (1) was proved in \cite{rezn2008,rezn2022-2008}. 

\begin{theorem}\label{ttg2}
Let $G$ be a regular paratopological group. 
If $G$ is $\lrd c\om$-Baire, then $G$ is a topological group. 
\end{theorem}
\begin{proof}
Let $O\in \cN_e$. Take  $\sqnn U\subset\cN_e$ such that $\cl{U_0^3}\subset O$ and $\cl{U_{n+1}^3}\subset U_n$ for $\nom$. Let us denote $P=\cN(O)$, $P_n=\cN(U_n)$ for $\nom$, $\cp= \sqnn P$, and $S=\ppl^c_\om(\cp)$. 
Statement \ref{s:tg:1} implies $\cl{\cN(P_0)\circ \cN(P_0)}\subset P$ and $\cl{\cN(P_{n+1})\circ \cN(P_{p+1})}\subset \cN(P_n)$ for $\nom$. Then $S\subset P$. Since $S$ is closed and $G$ is $\lrd c\om$-Baire, then $W\times W\subset P$ for some open nonempty $W\subset G$. Then $V=V^{-1}\subset O$ for $V=WW^{-1}$.
\end{proof}

\begin{proposition}[Theorem 4 \cite{Korovin1992}]\label{ptg3}
Let $X$ be a space such that $X^\om$ is a pseudocompact space. Let $H$ be an Abelian group such that $|X|\leq |H|=|H|^\om$.
Then there exists a semitopological group $G$ such that $G^\om$ is a pseudocompact space, $G$ maps continuously onto $X$, and there exists a homomorphism $H\to G$.
\end{proposition}

Recall that a group is called \term{Boolean} if it satisfies the identity $g^2=e$. In a Boolean group, taking the inverse element is the identity operation. Let $X$ be a compact space with uncountable Souslin number and $H$ a Boolean group such that $|X|\leq |H|=|H|^\om$. Then from Proposition \ref{ptg3} we obtain the following corollary.

\begin{cor}[Korovin] \label{ctg1}
There exist a pseudocompact Boolean quasitolpological group which is not a topological group.
\end{cor}

\section{Examples and questions}\label{sec-qe}

In this section, we study how different the introduced classes of spaces are. This section supplement the corresponding section in \cite{rezn2022gbtg}.

The following diagram shows the relationship between the most interesting classes of spaces.

Any arrow from $A$ to $B$ means that any $A$-Baire space is a $B$-Baire space.
Besides, arrow $\begin{tikzcd} A \arrow[r] & B \end{tikzcd}$ means that the converse is not true,  
arrow $\begin{tikzcd} A \arrow[r,dashrightarrow] & B \end{tikzcd}$ means that it is not known whether or not the converse is true, and arrow $\begin{tikzcd} A \arrow[r,tail, dashed] & B \end{tikzcd}$ means that there is a counterexample for the converse assertion under additional set-theoretic asseumptions.
In this diagram, a $\gbm$-Baire space is exactly a Baire space.

The following diagram summarizes 
\cite[Proposition 14]{rezn2022gbtg},
 Proposition \ref{pdsp1}, and Theorem \ref{tgd1}. 
Counterexamples are given in \cite{rezn2022gbtg} and will be constructed below for the bottom line of the diagram.

{
\def\b#1{\gbm_{#1}}
\def\o#1#2{\god_{#1,#2}}
\def\a#1{\arrow[#1]}
\def\na#1{\arrow[#1,dashrightarrow]}
\def\nna#1{\arrow[#1,tail, dashed]}
\[
\begin{tikzcd}
&\o ol \na d & \o pl \a l \na d &  \b{f} \a{l} \nna{d} &
\\
\b{}&\D \a l & \D_h \a l& \D_s \a l & \mathrm{CDP}
\na{l}
\end{tikzcd}
\]
}

\begin{problem} \label{pqe000} 
Distinguish the following triples of classes:
\begin{itemize}
\item[\rm(1)]
$\god_{o,l}$-Baire, $\godt_{o,l}$-Baire, and $\D$-Baire spaces;
\item[\rm(2)]
$\god_{p,l}$-Baire, $\godt_{p,l}$-Baire, and $\D_h$-Baire spaces;
\item[\rm(3)]
$\gbm_{f}$-Baire, $\gbmt_{f}$-Baire, and $\D_s$-Baire spaces;
\item[\rm(4)]
\cdp/-Baire and $\D_s$-Baire spaces.
\end{itemize}
\end{problem}


From Assertion 11, Proposition 14, the diagram from Section 7 from \cite{rezn2022gbtg}, and Theorem \ref{tgd1}  the following statement follows.

\begin{statement}\label{aqe1}
Let $X$ be a regular space without isolated points.
\begin{itemize}
\item[\rm(1)]
If $X$ is $\gbm_f$-nonmeager, then $X$ contains points with a countable base.
\item[\rm(2)]
If $X$ is compact, then $X$ is $\D_h$-Baire.
\item[\rm(3)]
If $X$ is pseudocompact, then $X$ is $\D$-Baire.
%
\end{itemize}
\end{statement}

Below we give examples distinguishing between the classes of spaces listed above.
In the examples below, \bitem B,A. $X$ means that $X$ is an $A$-Baire space that is not $B$-nonmeager.

\begin{example}\label{eqe6}
\bitem \gbm_{f},\mathrm{CDP}. $\R^{\om_1}$ (assuming MA+$\lnot$CH).
\edemo
Assertion \ref{aqe1} (1) and Proposition \ref{pcdp2}.
\endedemo
\end{example}


We denote by $\B$ the discrete two-point space $\sset{0,1}$. The base of the topology in $\B^C$ is formed by sets of the form
\[
W(A,B,C)\eqdef \set{\sq{x_\al}{\al\in C}\in \B^C: x_\al=0 \text{ for }\al\in A\text{ and }x_\be=1 \text{ for }\be\in B }
\]
for finite disjoint $A,B\subset C$.

\begin{example}\label{eqe7}
\bitem \D_s, \D_h. $\B^{2^\om}$.
\edemo
We denote $X=\B^{2^\om}$.
The space $X$ is $\D_{h}$-Baire (Statement \ref{aqe1} (2)). Let us show that $X$ is not $\D_s$-nonmeager. There is a cover $\set{F_\al:\al<2^\om}$ of the space $X$ of size $2^\om$, consisting of nowhere dense closed sets. We fix $\al_x<2^\om$ for $x\in X$ such that $x\in F_{\al_x}$.
We set $\til x=\sq{\til x_\al}{\al\in 2^\om}$, where $\til x_\al=x_\al$ for $\al\neq \al_x$ and $ \til x_{\al_x}\neq x_{\al_x}$. We define an almost neighborhood $P$ of the diagonal $X\times X$,  by setting $P(x)=X\setminus\sset{\til x}$.
Let us show that $\sorestr P V\notin \s_s(V)$ for any open nonempty $V\subset X$. Let $U=W(A,B,2^\om)\subset V$ for some finite disjoint $A,B\subset 2^\om$. We set
\[
W=V\setminus \bigcup_{\al\in A\cup B} F_\al.
\]
Then $W$ is dense and open in $U$ and $\til x\in U$ for $x\in W$. Consequently
\[
M=\set{x\in U: U\subset P(x)} \cap W = \es
\]
and the set $M$ is nowhere dense in $U$.

Note that the compact space $\B^{2^\om}$ is a $\gbm_{k}$-Baire space (see \cite{rezn2022gbtg} for the definition of $\gbm_{k}$-Baire spaces).
\endedemo
\end{example}

\begin{example}\label{eqe8}
\bitem \D_h, \D.
$G$. Let $G$ be a pseudocompact Boolean quasitolpological group that is not a topological group (Corollary \ref{ctg1}).
\edemo
The space $G$ is a $\D$-Baire space (Statement \ref{aqe1} (3)).  Theorem \ref{ttg1} implies that $G$ is not a $\D_h$-nonmeager space.

Note that the pseudocompact space $G$ is a $\gbm_{o}$-Baire space (see \cite{rezn2022gbtg} for the definition of $\gbm_{o}$-Baire spaces).
\endedemo
\end{example}

Let us denote by $\Sl$ the  Sorgenfrey line, the line $\R$ with the topology whose base is formed by the half-open intervals $[a,b)$. 

\begin{example}\label{eqe9}
\bitem \D, \gbm. $\Sl$.
\edemo
The  Sorgenfrey line $\Sl$ is a regular Baire paratopological group which is not a topological group.
Theorem \ref{ttg1} implies that $\Sl$ is not a $\D$-nonmeager space. 
There are also immediate simple arguments that show that $\Sl$ is not $\D$-nonmeager.
Let $P=\set{(x,y)\in\Sl: y\geq x}$. Then $P$ is a closed semi-neighborhood of the diagonal and $U\times U\not \subset P$ for any nonempty open $U\subset \Sl$.
Note that $\Sl$ is not a $\lrd c{\la}$-nonmeager space for any ordinal $\la$.
\endedemo
\end{example}

\begin{example}\label{eqe10}
Let $T$ be a countable space whose nonempty open sets are complements to finite sets. Then $T$ is $\D_v$-Baire meager $T_1$ space.
\edemo
Obviously, $T$ is a meager $T_1$ space.
Since any open nonempty set in $T$ is dense in $T$, then $\clv P=T\times T$ for any $P\in\sod T$. Therefore $T$ is $\D_v$-Baire.
\endedemo
\end{example}

\begin{problem} \label{e:qe2:1}
Let $X$ is a regular $\D_a$-Baire space. Is it true that $X$ contains an isolated point?
\end{problem}

A space $X$ is called \term{weakly pseudocompact} if there exists a compact Hausdorff extension $bX$ of the space $X$ in which the space $X$ is $G_\de$-dense, i.e., $X$ intersects any nonempty $G_\de$ subset of $bX$ \cite{arh-rezn2005}. It is clear that a product of weakly pseudocompact spaces is weakly pseudocompact; in particular, a product of pseudocompact spaces is weakly pseudocompact.

\begin{problem} \label{pqe5}
Let $X$ be a weakly pseudocompact space (or a product of pseudocompact spaces). 
Is $X$ a $\D$-Baire space?
\end{problem}


\begin{problem} \label{pqe5+1}
Let $G$ be a 
 paratopological group which is a
weakly pseudocompact space (a product of pseudocompact spaces). Is it true that $G$ is a topological group?
\end{problem}

\begin{problem} \label{pqe6}
Let $X$ and $Y$ be (completely) regular countably compact spaces. To which of the following classes does the product $X\times Y$ belong: $\D_h$-Baire spaces, $\D$-Baire spaces?
\end{problem}

\begin{problem} \label{pqe7}
Are \cdp/-Baire spaces $\god_{o,l}$-Baire?
\end{problem}

\begin{problem} \label{pqe8}
Is it possible to construct, without additional set-theoretic assumptions, an example of a \cdp/-space without a dense metrizable Baire subspace?
\end{problem}

\begin{problem} \label{pqe9}
Let $1<\ga<\om$.
Are there $\rd \ga$-Baire and $\lrd c\ga$-Baire spaces that are not $\D$-Baire?
\end{problem}

\begin{problem} \label{pqe9+1}
Let $X$ be a regular $\lrd c\om$-Baire space. Which of the following conditions are true?
\begin{itemize}
\item[\rm(1)]
$X$ is Baire.
\item[\rm(2)]
If $X$ is metrizable then $X$ is Baire.
\item[\rm(3)]
$X$ is $\D$-Baire.
\end{itemize}
\end{problem}

Note that the rational numbers $\Q$ are not $\lrd c\om$-Baire. There are paratopological groups that are not topological groups and are homeomorphic to $\Q$. For example, countable dense subgroups of the Sorgenfrey line $\Sl$. It follows from Theorem \ref{ttg2} that $\Q$ are not $\lrd c\om$-Baire.

\begin{problem} \label{pqe10}
Distinguish between the classes of spaces that are introduced in Section \ref{sec-seminghbr2}.
\end{problem}

The author thanks the referee for useful comments.

\bibliographystyle{elsarticle-num}
\bibliography{gbd}
\end{document}